\documentclass[11pt]{amsart}
\usepackage{amsmath,amsfonts,amssymb,graphicx,amsthm,hyperref}
\usepackage{verbatim,bbm,latexsym,indentfirst}
\usepackage{enumerate}
\usepackage{epsfig,color}
\usepackage{epstopdf}

\newcommand{\geom}{\mathrm{geom}}

\def\Lra{ \Longrightarrow }
\def\mgf{moment generating function }  
\def\geo{geometric }
\def\sf{scale function } 
\def\bfun{barrier influence function }
\def\claim{C} 
\def\mD{\mathcal D}
\def\mC{{\mathcal C}}
\def\T{\widetilde}

 \def\mCB F{${\mathcal CB F}$}

\def\Rui{\Psi} \def\sRui{\bar{\Rui}} 
 \def\Rui{\Psi}

\newcommand{\be}{\begin{eqnarray}}
\newcommand{\ee}{\end{eqnarray}}
\newcommand{\ba}{\begin{array}}
\newcommand{\ea}{\end{array}}
\newcommand{\baa}{\left[\begin{array}}
\newcommand{\eaa}{\end{array}\right]}

\def\BEN{\begin{enumerate}}  \def\BI{\begin{itemize}}
\def\EEN{\end{enumerate}}   \def\EI{\end{itemize}}

\def\beq{\begin{equation}} \def\eeq{\end{equation}}
\def\bea{\begin{eqnarray*}}
\def\eea{\end{eqnarray*}}
\def\le{\left} \def\ri{\right}
\def\bar{\overline}

\def\I{\infty}  \def\a{\alpha}

  \def\th{\theta}
\def\tw{w} \def\Tt{\T \t_{-1}^-}

\def\t{\tau}   \def\f{\phi}   \def\vf{\varphi}
    \def\qu{\quad} \def\D{\Delta}
\def\F{\Phi}    
  
\def\td{\t^-_{-1}}

\def\la{\label} \def\fr{\frac} \def\im{\item}

\newtheorem{Thm}{Theorem}

\def\beT{\begin{Thm}}
\def\eeT{\end{Thm}}

\newtheorem{Qu}{Problem}
\def\beQ{\begin{Qu}}
\def\eeQ{\end{Qu}}

\newtheorem{Lem}{Lemma}
\newtheorem{Exa}[Lem]{Example}

\newtheorem{Cor}[Lem]{Corollary}
\newtheorem{Def}[Lem]{Definition}
\newtheorem{Rem}[Lem]{Remark}
\newtheorem{Exe}[Lem]{Exercice}
\def\beXe{\begin{Exe}} \def\eeXe{\end{Exe}}
\def\eeD{\end{Def}} \def\beD{\begin{Def}}
\def\beXa{\begin{Exa}} \def\eeXa{\end{Exa}}
\def\beR{\begin{Rem}} \def\eeR{\end{Rem}}
\def\beL{\begin{Lem}} \def\eeL{\end{Lem}}
\newtheorem{Pro}[Lem]{Proposition}
\def\beP{\begin{Pro}} \def\eeP{\end{Pro}}
\def\beC{\begin{Cor}}
\def\eeC{\end{Cor}}
\def\bc{\begin{cases}}
\def\ec{\end{cases}}
\def\BEN{\begin{enumerate}} \def\EEN{\end{enumerate}}

\def\ssec{\subsection} \def\sec{\section}

  \def\CL{Cram\'er-Lundberg }
  \def\PK{Pollaczek-Khinchine }

  \def\pgf{probability generating function }
  
    \def\P{P}
    \def\E{E}
    \def\Z{{\mathbb Z}} 
    \def\N{{\mathbb N}}

\def\gf{generating function }

\def\sn{spectrally negative }

\def\kil{\mathcal E}

\def\qu{\quad}   
\def\lev{L\'evy }     
\def\deF{de Finetti }  

 \def\fp{first passage }

    \def\Eb{\E^{b]}} \def\Enull{\E^{0]}}   
    \def\Ezb{\E^{[0,b]}}

     \def\deF{de Finetti }

 \def\pg{\T{p}}

\begin{document}
\title{First passage problems for  upwards skip-free random walks via the $\F,W,Z$ paradigm}

\author{Florin Avram}
\address{Laboratoire de Math\'ematiques Appliqu\'ees, Universit\'e de Pau,  France}
\email{florin.avram@univ-Pau.fr}
\author{Matija Vidmar}
\address{Department of Mathematics, University of Ljubljana, Slovenia}
\address{Institute for Mathematics, Physics and Mechanics, Ljubljana, Slovenia}
\email{matija.vidmar@fmf.uni-lj.si}

\thanks{MV acknowledges financial support from the Slovenian Research Agency (research core funding No. P1-0222).}

\date{\today}

\maketitle
\begin{abstract}
We develop the  theory of the $W$ and $Z$ scale functions for right-continuous (upwards skip-free) discrete-time discrete-space random walks, along the lines of the analogous theory for spectrally negative L\'evy processes. Notably, we introduce for the first time in this context the one and two-parameter scale functions $Z$, which appear for example in the joint problem of deficit at ruin and time of ruin, and in problems concerning the walk reflected at an upper barrier. Comparisons are made between the various theories of scale functions as one makes time and/or space continuous. The theory is shown to be  fruitful by providing a convenient unified framework for studying dividends-capital injection problems under various objectives, for the so-called compound binomial risk model of actuarial science.
\end{abstract}
\vspace{0.25cm}
\noindent {\bf Key words}: skip-free Markovian jump processes; random walks; scale functions; martingales; compound binomial risk model;  dividends; capital injections.

\vspace{0.25cm}
\noindent \textit{2010 Mathematics Subject Classification}: Primary: 60G50; Secondary: 91B30.
\tableofcontents

\sec{Introduction}
First passage theory for random walks is  a classic topic, excellently treated for example in the textbooks \cite{Spitzer,Fel,takacs1977combinatorial,Bor}, and this  includes the upwards skip-free compound binomial model of the actuarial literature. However, in light of recent developments in the parallel continuous-time theory of spectrally negative/upwards skip-free L\'evy and Markov additive processes --- see for example \cite{AKP,Iva,IP,vidmar2013fluctuation,AIZ,AGV} --- it seems worthwhile to revisit this topic.

Indeed, while it is well-known that optimization problems in the discrete setup (which is in many ways more natural than the continuous one) may be tackled  numerically by dynamic programming algorithms, it is less known that when restricting to the skip-free case, the solutions of a great variety of \fp problems may be parsimoniously expressed in terms of two families of scale functions, 
just like in the continuous-time \lev case.

Recall that in the \lev case the scale functions $W^{(q)}$ and $Z^{(q)}$ have been known since \cite{Suprun} and \cite{AKP}, and that these functions intervene in important optimization problems. For example, $W^{(q)}$ provides the value function of the classic de Finetti problem of optimizing expected dividends until ruin with discount factor $q$ \cite{APP}, and $Z^{(q)}(\cdot,\th)$  intervenes for instance in the moment generating function (as function of $\theta$) of the capital injections \cite{IP} and in the combined dividend payout-capital injections problem for a doubly reflected process \cite{APP,AIjoint}. These are just two examples from an ever increasing list of problems \cite{Pispot,Kyp,AGV}, which can be now tackled by simple lookup in the list and using off-shelf packages computing the functions $W$ and $Z$ \cite{Iva}.

It was expected that the \fp theory developed in the world of spectrally negative L\'evy processes, which we call the $\F,W,Z$ paradigm, should have parallels for other classes of spectrally negative/skip-free Markov processes. In particular, the three cases listed below, being precisely the processes with stationary independent increments that exhibit non-random overshoots \cite{vidmar2015overshoots} (modulo trivial processes with monotone paths), were expected to be very similar:
\begin{enumerate}[(i)]
\item\label{type:i} (discrete-time, discrete-space) right-continuous (i.e. skip-free to the right) random walks, also known in insurance as the compound binomial model;
\item\label{type:ii} (continuous-time, discrete-space) compound Poisson processes that live on a lattice $h\mathbb{Z}$, $h\in (0,\infty)$, jumping up only by $h$ (what were called upwards skip-free L\'evy chains in \cite{vidmar2013fluctuation});
\item\label{type:iii} (continuous-time, continuous-space) spectrally negative L\'evy processes.
\end{enumerate}

However, important steps were missing for the fully discrete setup. Notably, the second scale function $Z_v(\cdot,\tw)$ was absent from the previous literature, and we provide below for the first time its generating function ($z$-transform) \eqref{eq:Z-generating-fnct}.

A second contribution of our paper is spelling out the connections between the three types of  first passage problems listed above. In particular, we provide in Appendix~\ref{summary-table} a concise table featuring side-by-side some of the salient features of the $\F,W,Z$ theory for the three types of process \eqref{type:i}-\eqref{type:ii}-\eqref{type:iii} delineated above. It may serve as an inexhaustive summary and a quick reference; for the complete exposition, the main body of the text must be consulted.

A third contribution is showing the convenience of using  the $\F,W,Z$ theory for solving dividends-capital injections problems -- see Sections~\ref{section:application} and~\ref{sec:examples}.


Now, the doubly discrete (in time and space) random walk risk model is defined by \cite{gerber1988mathematical, shiu1989probability}:
$$ X_n =X_0 + c n- \sum_{i=1}^n \claim_i, \quad n\in \N_0,$$
where $X_0$, taking values in $\Z$, is the initial capital, $c \in \N$ is the premium rate and the $\claim_i $, $i\in \mathbb{N}$, take values in $\mathbb{N}_0$ and are independent, identically distributed random variables with probability mass function $p_k=P(\claim_1=k)$ for $k\in \mathbb{N}_0$. One advantage of the discrete setup over the more popular continuous
time models is the possibility to replace the Wiener-Hopf factorization by the conceptually simpler factorization of  Laurent series (see for example \cite{banderier2002basic} and \cite{xin2004ring}); another advantage is that one has access to Panjer recursions for computing compound distributions.


The results simplify considerably for the upwards skip-free compound binomial model obtained when $c=1$ (\cite{quine,brown,marchal2001combinatorial} \cite[\emph{passim}]{Spitzer} \cite[Section~4.1]{PP} among others):
 \beq \label{discmod} X_n =X_0 + n -  \sum_{i=1}^n \claim_i,\quad n\in \N_0,
\eeq
that we now consider as having been fixed and to which we specialize all discussion henceforth. We insist throughout that $p_0>0$. 

Notation-wise, we let $$\pg( z):=E z^{\claim_1}=\sum_{k=0}^\I p_k z^k,\quad z\in (0,1],$$
 denote the probability generating function of the claims. Then, for $n\in \mathbb{N}$, (in the obvious notation) $E z^{\sum_{i=1}^n \claim_i}=\big[\tilde p(z)\big]^n=(p_0 + (1-p_0) \pg_{\claim|\claim \geq 1}( z))^n$, which makes it manifest that $\sum_{i=1}^n \claim_i$, the total claims arising from $n$ time periods, has a compound binomial distribution, explaining the name compound binomial model: at each instant in discrete time, a positive claim either occurs or not, with probability $1-p_0$ and $p_0$, respectively, independently of the sizes of the positive claims.

 \beR  By the independence of the claims, we may also write, for $n\in \mathbb{N}_0$:
$$E
 \le[z^{\sum_{i=1}^n (\claim_i-1)}\ri]= \left(\fr{\pg( z)}{z}\right)^n \qu
  \Lra \qu \sum_{m=0}^\I v^m E
 \le[z^{\sum_{i=1}^m (\claim_i-1)}\ri]=\fr 1{1-v \pg( z)/z},\quad v\in \left(0,\frac{z}{\pg(z)}\right).$$
The last expression, called  the ``unrestricted generating function" in \cite[Eq.~(8)]{banderier2002basic}, identifies already potential singularities as the roots of the Lundberg equation \cite{Lun} $\fr {\pg( z)}{z}=v^{-1}$. The smallest (positive) root of this equation plays a central role in our story --- see next section. \eeR

Next, we will denote by
\begin{equation}
\tau_b^-=\inf\{t\geq 0: X_t\leq b\}\text{ and } \tau_b^+=\inf\{t\geq 0: X_t\geq b\},
\la{tau}
\end{equation}
respectively, the first passage times  below and above a level $b$  (with $\inf\emptyset=\infty$).
\beR
Note this differs slightly from the usual definition of these quantities for a spectrally negative L\'evy process, say $U$. There one replaces $t\geq 0$ by $t>0$ and $\leq b$ ($\geq b$) by $<b$ ($>b$); and, of course, $X$ by $U$. When considering $\tau^{\pm}_b$ for a spectrally negative L\'evy process $U$, we shall mean these quantities with the latter replacements having been effected.
\eeR

 Lastly, for convenience, we assume given a family of measures $(P_x)_{x\in \mathbb{Z}}$ with corresponding expectation operators $(E_x)_{x\in \mathbb{Z}}$, for which: (i) $P_x(X_0=x)=1$ for all $x\in \Z$; and (ii) the $C_i$, $i\in \mathbb{N}$, have the same law under all the $P_x$, $x\in \Z$, as they do under $P=P_0$.

\beR\label{remark:connection:1}
The discrete-time discrete-space compound binomial model is embedded into continuous time via subordination (time-change) by an independent homogeneous Poisson process $N$. In precise terms, allowing also a scaling of space, we have the following correspondence between the right-continuous random walk $X$ of \eqref{discmod} and the upwards skip-free L\'evy chain of \cite[Sec.~2]{vidmar2013fluctuation} that we will here denote by $Y$:
$$ X\quad \rightsquigarrow\quad Y:\quad  Y_t:=hX_{N_t},\quad t\in[0,\infty),$$
where $h\in (0,\infty)$ is space scaling. In particular, denoting the intensity of $N$ by $\gamma$, the L\'evy measure $\lambda$ of $Y$ is given by $\lambda=\gamma\sum_{i\in \mathbb{Z}\backslash \{1\}}p_i\delta_{h(1-i)}$; and if we denote the Laplace exponent of $Y$ by $\psi$ (so $\psi(\beta)=\frac{\log E[e^{\beta Y_t}]}{t}$ for $\beta\in [0,\infty)$), then $\psi(\beta)=\gamma[e^{\beta h}\pg(e^{-\beta h})-1]$. Note that the mass of the L\'evy measure $\lambda$ is $\gamma(1-p_1)$, which may be strictly less than $\gamma$.
\eeR

\beR
In the following, when the $\tau_b^\pm$ appear in the context of the upwards skip-free L\'evy chain $Y$, they are to be interpreted in the sense of \eqref{tau} with $Y$ replacing $X$.
\eeR

Here is now a brief guide to the contents. In Sections~\ref{sec:smooth},~\ref{sec:non-smooth} and~\ref{sec:smooth-two-sided}, we review, respectively (with $v$ indicating discounting):
 \BEN \im  the smooth one-sided first passage problem, which introduces the Lundberg root $\vf_v$ (analogue of $\F(q)$ from the \lev theory);
\im the non-smooth one-sided first passage problem, which involves the ruin and survival probabilities $\Rui_v$, $\sRui_v$;
 \im the smooth two-sided first passage problem, where the fundamental scale function $W_v$ first appears.
 \EEN

We turn then to new material in Section~\ref{sec:deficit-at-ruin-Z}, computing the generating function ($z$-transform) of the second hero of \fp theory: the $Z_v(\cdot,\tw)$ scale function. This is introduced via the problem of deficit at ruin: we provide the analogue \eqref{eq:Z-via-deficit-at-ruin} of the following two-sided exit identity for a spectrally negative L\'evy processes (in standard notation):
$$E_x [e^{-q \t^{-}_0 +\theta X(\t^-_0)} ; \t^-_{0} <\t_b^+]=Z^{(q)}(x,\theta)-\fr{W^{(q)}(x)}{W^{(q)}(b)}Z^{(q)}(b,\theta), \la{bea}$$
 with its beautiful probabilistic interpretation \cite[Cor.~3]{IP}. We also determine the analogue \eqref{eq:levy-identity} of the formula \cite[Eq.~(8.9)]{Kyp} (again for a spectrally negative L\'evy process, in standard notation) $$
E_x [e^{-q \t^{-}_0} ; \t^-_{0} <\infty]=
 Z^{(q)}(x) - \fr q{\F(q)} W^{(q)} (x),\quad q>0,$$ which is interesting, for example, since it reveals that the two protagonists of the ``reflected" and ``absorbed" smooth passage problems, $Z^{(q)}$ and  $W^{(q)}$, have the same asymptotics at $\infty$, up to a constant.  A distinguishing element of the scale functions $W_v$ and $Z_v(\cdot,w)$, in the present context, are explicit recursions available for their computation: see \eqref{har} and \eqref{eq:Z-recurrence}, respectively. Section~\ref{section:application} discusses some important applications, like the  \deF dividends optimization problem, and the optimization of dividends for the doubly reflected process. These are complemented by illustrative numerical examples in Section \ref{sec:examples}. Finally, note that while our motivation for this investigation comes chiefly from risk models in the insurance context, the results presented are general and hence more widely applicable. 

\section{Smooth one-sided first passage problem: the Lundberg equation}\label{sec:smooth}

The first  key observation is that for the first passage upwards, the stationary independent increments and skip-free properties imply a multiplicative structure; thus, for integer $x \leq b$, and for $v\in (0,1]$, we have
\begin{equation}\label{eq:smooth-one-sided-identity} E_x \le[v^{\t_b^+}; \t_b^+<\infty\ri]=\vf_v^{b-x}, \end{equation}
where
$$\vf_v:=E \le[v^{\t_1^+};\t_1^+<\infty\ri]= \sum_{k=1}^\I v^k P[\t_1^+=k]  \in (0,v].$$

Conditioning at time $1$, we obtain $$\vf_v =
v\, E \le[E_{ 1 -\claim_1}[ v^{\t_1^+};\t_1^+<\infty]\ri]=v\, \sum_{k=0}^\I  p_k  \vf_v^{k}=v
 \pg(\vf_v),$$ which   reveals that $\vf_v$ appearing in \eqref{eq:smooth-one-sided-identity} satisfies the   Lundberg equation \cite[Eq. ~(3.3)]{cheng2000discounted}, \cite[Eq.~(6.8)]{gerber2010elementary}: \beq \la{disL}\fr{\vf_v}{\pg(\vf_v)}=v. \eeq
Alternatively, this relation may be derived by looking for exponential martingales of the form $(v^t \xi^{-X_t})_{t\in \N_0}$, for fixed $v$, and $\xi$ from $(0,1]$: $(v^t\xi^{-X_t})_{t\in \N_0}$ is a martingale iff $\fr{\xi}{\pg(\xi)}=v$; and then applying optional sampling.
\beR
The function $(0,1]\ni \xi\mapsto \pg(\xi)/\xi=E \xi^{\claim_1-1}$ is strictly convex, equal to $1$ at $1$, and tending to $\infty$ at $0$. It follows that the equation (in $\xi\in (0,1]$) $\frac{\pg(\xi)}{\xi}=v^{-1}$ has a unique solution $\vf_v\in (0,1)$, when $v<1$ (furthermore, in this case, $\vf_v <v$), whereas in the case $v=1$, this equation has 
one or 
two solutions (one of which is always $1$), according as to whether $E \claim_1\leq 1$ or $E \claim_1>1$.
In the latter case $X$ drifts to $-\infty$, and $\vf_1 \in (0,1)$ is the smallest solution to $\xi=\pg(\xi)$ (in $\xi\in (0,1]$). Altogether,  this defines a continuous strictly increasing function $\vf: (0,1]\to (0,\vf_1 ]$.
\eeR
\beR
If for $q\in [0,\infty)$, we let $\Phi(q)$ be the largest zero of $\psi-q$, then  we see from Remark~\ref{remark:connection:1} that $\vf_v=e^{-h\Phi(\gamma(v^{-1}-1))}$ for all $v\in (0,1]$.
\eeR
\beR Note that \eqref{disL} identifies $\t_1^+$ as a Lagrangian type distribution \cite{consul2006lagrangian}. Indeed the distribution of $\t_1^+$ may be obtained using the
Lagrange inversion formula
$$\vf_v=\sum_{n=1}^\I \fr{v^n}{n!}\left[\left(\frac{d}{dw}\right)^{n-1} \pg(w)^n\right]_{w=0} =\sum_{n=1}^\infty \frac{v^n}{n}p^{n*}(n-1), \label{eq:LI}
$$
where for $n\in \mathbb{N}$, $p^{n*}$ is the $n$-fold convolution of the distribution $p$ with itself.
More generally, for $b\in \mathbb{N}$,
$$\vf_v^b=b\sum_{n=b}^\infty \frac{v^n}{n}p^{n*}(n-b),$$
 yielding Kemperman's formula \cite{kemperman1961passage} for the distribution of $\t_b^+$:
$$ \la{K} P[\t_b^+=n]=\frac{b}{n}p^{n*}(n-b)=\fr{b}{n} P[X_n=b],\quad n\in \N_{\geq b}.
$$
\eeR

\section{Non-smooth one-sided first passage problem: ruin and survival probabilities; the Lundberg  recurrence}\label{sec:non-smooth}
For initial capital $x\in \Z$, the finite time and eventual
ruin probabilities are defined by: $$ \Rui(n;x):=P_x [\t_{-1}^- \leq n]\text{ for } n\in \N_0,\qu \Rui(x):=\lim_{n\to\I}\Rui(n;x)=P_x [\t_{-1}^- < \I];$$
similarly we introduce the finite time and  perpetual survival probabilities: $$ \sRui(n;x):=P_x [\t_{-1}^- > n]\text{ for } n\in \N_0, \qu  \sRui(x):=\lim_{n \to \I} \sRui(n;x)=P_x [\t_{-1}^- =\infty].$$ Of course $\Rui(n;x)+\sRui(n;x)=1$,  $\Rui(x)+\sRui(x)=1$, and one has the recursions, valid for all integer $x  \geq 0$, $ n \geq 1$:
\begin{equation} \la{dreci1}
\sRui(n;x)=\sum_{i=0}^{x + 1} p_i  \sRui(n-1;x+ 1 -i),\; \sRui(0;x)=1,
\end{equation}
\begin{equation}\la{dreci2}
 \Rui(n;x)=\sum_{i=0}^{x + 1} p_i  \Rui(n-1;x+ 1 -i) + \sum_{k=x+2}^\infty p_k,\; \Rui(0;x)=0.
\end{equation}
These two recurrences may, for a sequence of functions $f_n:\Z\to [0,1]$, standing in lieu of $\Rui(n;\cdot)$, $\sRui(n;\cdot)$, be written symbolically as 
$$f_n= K \pg(K^{-1}) f_{n-1}\text{ on }\mathbb{N}_0,$$ which passes to the limit (as $n\to\infty$)
\begin{equation}\label{eq:Lundberg}
 f= K \pg(K^{-1}) f \text{ on }\mathbb{N}_0,
\end{equation}
 where $K$ is the translation operator, $Kg(x) :=g(x+1)$, and $f(x):=\lim_{n\to\infty}f_n(x) $. This limiting recurrence (satisfied by the eventual ruin and perpetual survival
probabilities $\Rui$ and $\sRui$) 
may be called the
  ``Lundberg
recurrence". It constitutes a linear difference equation for $f$, whose characteristic equation  is (in $x\ne 0$) $1=x\pg(1/x)$. The latter is (formally) just the
Lundberg equation \eqref{disL} with $v=1$ upon substituting $x^{-1}$ for $\vf_1$. When the distribution $p$ has a finite support, then from the theory of finite order linear difference equations with constant coefficients, this implies that $f$, in particular the ultimate ruin and  perpetual survival probabilities, may be expressed as combinations of powers of the roots of the characteristic equation (in $x\ne 0$)
\begin{equation}\label{eq:characteristic-general}
1=x\pg(1/x).
\end{equation}

Classical ruin theory proceeds by computing double (generating function) transforms, briefly reviewed in Appendix~\ref{s:Wil}.
For example, one useful result, similar to the \PK formula for the \CL model, is \cite[Eq.~(3.5)]{willmot1993ruin}\footnotesize
\begin{equation}  \label{eq:perp-survival-transform}\T{\sRui}(z):=\sum_{x=0}^\I   z^x  \sRui(x)=\fr{(1-E [\claim_1])\lor 0}{\pg(z)-z},\qu z\in (0,1).
\end{equation}\normalsize
Another is
\footnotesize
\begin{equation}\label{eq:ruin-transform}  \T{\Rui}_v(z):=\sum_{x=0}^\I \sum_{n=0}^\I  z^x v^n \Rui(n;x)=\fr{1 }{{z}-v {\pg(z)}}\le(\fr{v(z-\pg(z))}{(1-v)(1-z)}+\fr{\vf_v}{1-\vf_v}\ri), \quad v, z \in (0,1),z\ne\vf_v.\end{equation}
\normalsize

We will follow next an alternate approach, which focuses on the two-sided exit problem from an interval.

\section{Smooth two-sided first passage problem: the $W$ scale functions}\label{sec:smooth-two-sided}
In the context of L\'evy processes, the $W^{(q)}$ scale function is often defined first for $q=0$, in the case when the underlying process drifts to $\infty$, by proportionality to the survival probability, and then in the remainder of the cases by an Esscher transform/approximation \cite[Sec.~VII.2]{Ber97}  \cite[Sec.~8.2]{Kyp} \cite[Sec. 4.2]{vidmar2013fluctuation}. 

In our setting of the right-continuous random walk $X$, we introduce, for $v\in (0,1]$, the discrete-time analogue $W_v$ of $W^{(q)}$, by setting $W_v(y):=(p_0E[v^{ \t_y^+};  \t_y^+ < \t_{-1}^-])^{-1}$ for $y\in \mathbb{N}_0$ and $W_v(y)=0$ for $y\in -\mathbb{N}$. The Markov property at the time $\tau_x^+$ and the skip-free property (yielding $X_{\tau_x^+}=x$ on $\{\tau_x^+<\infty\}$) then imply the ``gambler's winning" relation \cite{marchal2001combinatorial,gerber2006note}, for integer $x \leq N$, $0\leq N$:
\beq E_x [v^{ \t_N^+};  \t_N^+ < \t_{-1}^-]=\fr{ W_v(x)}{ W_v(N)}. \la{Wdef} \eeq
We call $W_v$ the $v$-scale function and we write simply $W$ for the $1$-scale function $W_1$. (The choice of the normalization $W_v(0)=1/p_0$ is somewhat arbitrary, though it is guided by obtaining the simplest possible form for the $z$-transform of $W_v$ (\eqref{Wdisc} below); by comparison to the $W$ scale function of \cite{vidmar2013fluctuation} (see Remark~\ref{remark:Ws-connection} below); and the simplicity of subsequent formulae in which $W_v$ features.)
\beR
We use the subscript notation $W_v$ for the scale functions of $X$, reserving the superscript version $W^{(q)}$ for the corresponding quantities from the L\'evy setting. When only $W$ appears, it will be clear from context which of the two is meant. We will adhere to a similar convention with respect to the scale functions $Z^{(q)}(\cdot, \theta)$, $Z^{(q)}:=Z^{(q)}(\cdot, 0)$ and (hence the notation) their discrete-time analogues $Z_v(\cdot,\tw)$, $Z_v$.
\eeR

Conditioning on the first jump, \eqref{Wdef} implies the harmonic recursion
\cite[Eq.~(3.1)]{marchal2001combinatorial}
\beq \la{har}  W_v(x)=v \sum_{y=-1}^{x} W_v(x-y) p_{y+1}, \quad x\in \mathbb{N}_0.\eeq

Taking $z$-transform yields \cite[Eq.~(3.2)]{marchal2001combinatorial}
\beq \T W_v(z):=\sum_{x=0}^\I z^x W_v(x)=\fr {1}{\pg(z)-\frac{z}{v}}, \;\quad z \in (0,\vf_v).\la{Wdisc} \eeq

Since the $z$-transform \eqref{Wdisc} of $W_v$ is known,
  the computation of the \sf $ W_v$  reduces  finally
  to Taylor coefficient extraction of \eqref{Wdisc} expanded in
   a power series.

\beR\label{remark:killing}
It is seen from \eqref{har}, or directly from \eqref{Wdef}, that one has $W_v/v={^vW}$, where $^vW$ is the $1$-scale function of the process $X$ geometrically killed with probability $1-v$, i.e. of the process which has, \emph{ceteris paribus}, the sub-probability pmf $(vp_k)_{k\in \mathbb{N}_0}$ governing the sizes of the $\claim_n$, $n\in \mathbb{N}$. 
\eeR

\beR For $X$ embedded into continuous time as an upwards skip-free L\'evy chain, i.e. for the process $Y$ of Remark~\ref{remark:connection:1}, \eqref{har} and \eqref{Wdisc} become, respectively,   \cite[Eqs.~(4.10)~\&~(4.6)]{vidmar2013fluctuation}. This is seen through the identification $W^{(q)}(mh)=\frac{1}{\gamma h}W_{\frac{\gamma}{\gamma+q}}(m)$  for $m\in \mathbb{N}_0$, $q\in [0,\infty)$, where $W^{(q)}$ is the $q$-scale function of \cite{vidmar2013fluctuation}. Note also that the normalization $W_v(0)=p_0^{-1}$ is consistent with $W^{(q)}(0)=1/(h\lambda(\{h\}))=1/(\gamma hp_0)$ of \cite[Prop. (4.7)]{vidmar2013fluctuation}. On the other hand, in the spectrally negative case, there is no direct analogue of recursion \eqref{har}, though one can consider the heuristic relation (it is rigorous in the upwards skip-free case  \cite[Rem.~4.16]{vidmar2013fluctuation}) $(L-q)W^{(q)}=0$ on $(0,\infty)$ \cite[p.~136]{KKR}, $L$ being the infinitesimal generator of the underlying L\'evy process, to be a close relative. \eqref{Wdisc} has the Laplace transform equivalent \cite[Eq.~(8.8)]{KKR} that formally differs from \cite[Eq.~(4.6)]{vidmar2013fluctuation} only by the factor $(e^{\beta h}-1)/(\beta h)\to 1\text{ as }h\downarrow 0$ (with $\beta$ the argument of the Laplace transform).\label{remark:Ws-connection}
 \eeR

\beR \label{remark:alt-form-recursion-W}
An alternative form of  recursion \eqref{har} is \cite[Eq.~(4.13)]{vidmar2013fluctuation} $$W_v(n+1)=W_v(0)+\sum_{k=1}^{n+1}\frac{\frac{1}{v}-\sum_{l=0}^kp_l}{p_0}W_v(n+1-k),\quad n\in \mathbb{N}_0.$$ 
 In particular we see via induction that for each fixed integer $x$, the map $[1,\infty)\ni \xi\mapsto W_{1/\xi}(x)$ extends to a polynomial function defined on the whole of the complex plane.
\eeR

\beR
When $X$ drifts to $\infty$, i.e. when $E \claim_1 <1$, then with $v=1$, \eqref{Wdisc} coincides up to a multiplicative constant with the perpetual survival transform \eqref{eq:perp-survival-transform}.
We conclude that
$$\sRui(x)=(1- E(\claim_1))W(x).$$
\eeR

\beR\label{rmk:asymptotics}
It follows from \eqref{Wdisc} that $_vW(x)=\frac{\vf_v}{v}W_v(x)\vf_v^x$, where $_vW$ is the $1$-scale function of the Esscher transformed process in which $\claim_1$ has the geometrically tilted probability mass function $\mathbb{N}_0\ni k\mapsto \frac{v}{\vf_v}p_k\vf_v^k$. Hence by monotone convergence, $\lim_{x\to \infty}W_v(x)\vf_v^{x+1}=v\lim_{x\to\infty}{}_vW(x)=v\lim_{z\uparrow 1}\sum_{x=0}^\infty(1-z) z^x{}_vW(x)=\frac{v}{1-v\pg'(\vf_v-)}$, where we understand $1/0=\infty$. This confirms \cite[Prop.~4.8(i)]{vidmar2013fluctuation}. For a more detailed study of the behaviour of $W_1$ in the case when $\pg'(1-)=1$ and $\vf_1=1$, i.e. when $X$ oscillates, see \cite[Prop.~4.8(ii)]{vidmar2013fluctuation}.
\eeR

\beR
We note the  following interesting
observation of \cite{marchal2001combinatorial} that the scale function is essentially a determinant. For an arbitrary homogeneous Markov chain $(V_n)_{n\in \mathbb{N}_0}$ on a countable state space, let $(V_n')_{n\in \mathbb{N}_0}$ denote the chain killed outside a finite non-empty set $M$, and let $Q$ denote the corresponding restriction of the transition matrix to $M$. For $v\in (0,1)$, denote by $D_v$ the determinant of the  matrix $I-v Q$. Then the killed resolvent expresses as
$$\sum_{n=0}^\I P_i[V_n'=j] v^n=((I-v Q)^{-1})_{ij}=\fr{N_{ij}(v)}{D_v},\quad \{i,j\}\subset M,$$ where $N_{ij}(v)$ are the entries of the adjoint matrix $\mathrm{adj}(I-v Q)$ (see for example \cite[Cor. 2.2]{marchal2001combinatorial}). Restricting now to the upwards skip-free case (while \cite{marchal2001combinatorial} considers the downwards skip-free case), let, for $v\in (0,1]$, $D_v(N)$, $N \in\N$, denote the determinant corresponding (in the above sense) to the restriction  of $X$ to $\{0,1,2,...,N-1\}$,
and set $D_v(0):=1$. From \cite[Prop. 3.3]{marchal2001combinatorial},
$$ \la {twos} E_i[v^{\tau_N^+}, \tau_N^+<\tau_{-1}^-]=(p_0 v)^{N-i}\fr{D_v(i)}{D_v(N)},\quad \{i,N\}\subset \mathbb{N}_0,\; v\in (0,1).
$$
It follows that $W_v(i)=p_0^{-1}(p_0v)^{-i}D_v(i)$ for all $i\in \mathbb{N}_0$, $v\in (0,1]$.
\eeR

\beR
For $N\in \mathbb{N}$, the resolvent of the process $X$ killed on exiting $I_N:=\{0,\ldots,N-1\}$, denoted $X'$, is given by \cite[Prop.~3.2]{marchal2001combinatorial}
$$\sum_{n=0}^\I P_i[X_n'=j] v^n=v^{-1}\left(\frac{W_v(N-1-j) W_v(i)}{W_v(N)} -W_v(i-j-1)\right),\quad  \{i,j\}\subset I_N,\; v\in (0,1].$$ For the analogue of the latter in the spectrally negative case see e.g. \cite[Thm.~8.7]{Kyp}.
\eeR

We conclude this section with the important observation that

\beP\label{prop:scales:miscellaneous}
For each $v\in (0,1]$, $( v^{n\land \tau_{-1}^-} W_v(X_{n\land \tau_{-1}^-}))_{n\in \mathbb{N}_0}$ is a martingale under each $\P_x$, $x\in \mathbb{Z}$.
\eeP
\begin{proof}
This follows  from the harmonic recurrence \eqref{har}.
\end{proof}

\beR
The analogue of Proposition~\ref{prop:scales:miscellaneous} in the setting of upwards skip-free L\'evy chains are the martingales, for $q\in [0,\infty)$, $(e^{-q(t\land \tau_{-h}^-)}W^{(q)}(Y_{t\land \tau_{-h}^-}))_{t\in [0,\infty)}$ \cite[Cor.~4.17]{vidmar2013fluctuation}. In the case of a spectrally negative L\'evy process $U$,  $(e^{-q(t\land \tau_0^-)}W^{(q)}(U_{t\land \tau_0^-}))_{t\in [0,\infty)}$ is a local martingale with localizing sequence $(\tau_n^+)_{n\in \mathbb{N}}$ \cite[Ex.~8.12]{Kyp}. There are no issues with integrability in the discrete space case, because thanks to the skip-free property, $\P_x$-a.s. for any $x\in \mathbb{Z}$, by any deterministic time, the stopped process $X^{\tau_{-1}^-}$ is automatically bounded /and, for the upwards skip-free L\'evy chain $Y$, the further subordination by the independent homogeneous Poisson process $N$ does not ruin this/.
\eeR

\beC\label{corollary:aduivat}
For each $v\in (0,1]$ and integer $x\leq N$, $b\leq N$,
$$\E_x(W_v(X_{\tau_{b-1}^-})v^{\tau_{b-1}^-};\tau_{b-1}^-<\tau_N^+)=W_v(x)-\frac{W_v(x-b)}{W_v(N-b)}W_v(N).$$ In particular, $\E_x(W_v(X_{\tau_{b-1}^-})v^{\tau_{b-1}^-};\tau_{b-1}^-<\infty)=W_v(x)-W_v(x-b)\vf_v^b$.
\eeC
\begin{proof}
For any integer $x$, by optional sampling, the skip-free property and spatial homogeneity, $W_v(x)=E_x[W_v(X(\tau_{b-1}^-))v^{\tau_{b-1}^-};\tau_{b-1}^-<\tau_N^+]+E_x[W_v(X(\tau_N^+))v^{\tau_N^+};\tau_N^+<\tau_{b-1}^-]=E_x[W_v(X(\tau_{b-1}^-))v^{\tau_{b-1}^-};\tau_{b-1}^-<\tau_N^+]+W_v(N)E_x[v^{\tau_N^+};\tau_N^+<\tau_{b-1}^-]=E_x[W_v(X(\tau_{b-1}^-))v^{\tau_{b-1}^-};\tau_{b-1}^-<\tau_N^+]+W_v(N)E_{x-b}[v^{\tau_{N-b}^+};\tau_{N-b}^+<\tau_{-1}^-]$. The first identity then follows from \eqref{Wdef}. In particular, letting $N\uparrow \infty$ and using Remark~\ref{rmk:asymptotics}, we obtain the second identity (for instance first for $v<1$ and then taking the limit $v\uparrow 1$).
\end{proof}
\section{Problem of deficit at ruin with killing at an upper boundary: the $Z$ scale functions}\label{sec:deficit-at-ruin-Z}

Let $v\in (0,1]$, $w\in (0,1]$. For integer $x\leq b$, $b\geq 0$, by the Markov property at time $\tau_b^+$ and the skip-free property (yielding $X_{\tau_b^+}=b$ on $\{\tau_b^+<\infty\}$),
$$E_x [v^{ \t^-_{-1}} w^ {-X(\t^-_{-1})} ; \t^-_{-1} <\t_b^+]=E_x [v^{ \t^-_{-1}} w^ {-X(\t^-_{-1})} ; \t^-_{-1} <\infty]-E_x [v^{ \t^-_{-1}} w^ {-X(\t^-_{-1})} ; \t^+_b <\t_{-1}^-<\infty]$$
$$=E_x [v^{ \t^-_{-1}} w^ {-X(\t^-_{-1})} ; \t^-_{-1} <\infty]-E_x[v^{\tau_b^+};\tau_b^+<\tau_{-1}^-]E_b [v^{ \t^-_{-1}} w^ {-X(\t^-_{-1})} ; \t_{-1}^-<\infty].$$
Putting $\Rui_v(x,w):=E_x [v^{ \t^-_{-1}} w^ {-X(\t^-_{-1})} ; \t^-_{-1} <\I]$, we have then from the preceding and using \eqref{Wdef}, the neat identity
$E_x [v^{ \t^-_{-1}} w^ {-X(\t^-_{-1})} ; \t^-_{-1} <\t_b^+]=\Rui_v(x,w)-\frac{W_v(x)}{W_v(b)}\Rui_v(b,w)$. We introduce now, for some $\alpha_v(w)\in [0,\infty)$ that we shall specify in the sequel,
\begin{equation}\label{eq:Zdef}
Z_v(x,w):=\Rui_v(x,w)+\alpha_v(w) W_v(x),
\end{equation}
a slightly modified $\Rui_v(\cdot,w)$, which also satisfies the identity
\begin{equation}\label{eq:Z-via-deficit-at-ruin}
E_x [v^{ \t^-_{-1}} w^ {-X(\t^-_{-1})} ; \t^-_{-1} <\t_b^+]=Z_v(x,w)-\frac{W_v(x)}{W_v(b)}Z_v(b,w)
\end{equation}
(easy to check). The first motivation for preferring to use $Z_v(\cdot,w)$ with a suitable choice of $\alpha_v(w)$ instead
of $\Rui_v(\cdot,w)$ appears below in \eqref{eq:Z-generating-fnct}, and then in Section~\ref{section:application}; many other formulas where the analogue of $Z_v(\cdot,w)$ is preferable are known in the  literature on \sn \lev processes -- see for example \cite{IP,AGV}.

\beR \la{Zbc} Note that $Z_v(x,w)=\Rui_v(x,w)=w^{-x} \text{ for all integer } x \leq -1$.
\eeR
We compute now the $z$-transform of $Z$. Conditioning on the first jump, we obtain from \eqref{eq:Zdef} and the definition of $\Rui_v(\cdot,w)$, via \eqref{har}, the recurrence relation
\begin{equation}\label{eq:Z-recurrence}
Z_v(x,w)/v=\sum_{k=-1}^x p_{k+1}Z_v(x-k,w) + \sum_{k=x+1}^\I w^{k-x} p_{k+1},\quad x\in \mathbb{N}_0.
\end{equation}
Hence the generating function $\T Z_v(z,w):=\sum_{x=0}^\I z^x  Z_v(x,w)$ satisfies, for  $z\in (0,\vf_v)\backslash \{w\}$,
$$
\T Z_v(z,w)/v=p_0 \fr{\T Z_v(z,w)-Z_v(0,w)} z + \sum_{x=0}^\I z^x \sum_{k=0}^x Z_v(x-k,w) p_{k+1} + \sum_{x=0}^\I z^x \sum_{k = x+1}^\infty w^{k-x} p_{k+1}$$
$$=p_0 \fr{\T Z_v(z,w)-Z_v(0,w)} z + \sum_{k=0}^\I p_{k+1} z^k  \sum_{x=k}^\I z^{x-k} Z_v(x-k,w)  + \sum_{k=1}^\I p_{k+1} w^k  \sum_{x=0}^{ k-1}  \left(\fr z w\right)^{x}$$
$$ =p_0 \fr{\T Z_v(z,w)-Z_v(0,w)} z + \T Z_v(z,w) \sum_{k=0}^\I p_{k+1} z^k  + \sum_{k=1}^\I p_{k+1} w^k   \fr{1-(\fr z w)^{k}}{1-\fr z w}
$$
$$ =p_0 \fr{\T Z_v(z,w)-Z_v(0,w)} z + \T Z_v(z,w) \frac{\pg(z)-p_0}{z} +\fr{\frac{\pg(w)-p_0}{w}-\frac{\pg(z)-p_0}{z}}{1-\fr z w},
$$
i.e.,  in view of \eqref{Wdisc},
$$ \T Z_v(z,w)=
-p_0 (1-Z_v(0,w))\T W_v(z)+    \fr{z\T p (w)- w \T p(z)}{ (z -w)(\T p(z)-\frac{z}{v})}.$$

Recall now that in the \lev case, $Z^{(q)}(0,\theta)$ is chosen so as to ensure a ``smooth fit" \cite[Def.~5.8]{APP15} to the  boundary condition $e^{x \theta}$ for $x\in (-\infty,0)$. The analog in the discrete case  is to insist on
$Z_v(0,w)=1$, which we may do
by an appropriate  choice of $\alpha_v(w)$. Furthermore, this choice (that we assume henceforth) leads to the simple expression
\begin{equation}\label{eq:Z-generating-fnct}
\T Z_v(z,w)=\fr 1 {\T p(z)-\frac{z}{v}}
\fr{z\T p (w)- w \T p(z)}{ z -w},\quad z\in (0,\vf_v),\; v\in (0,1],\; w\in (0,1]
\end{equation}
(where the quotient must be understood in the limiting sense when $z=w$).

Extracting the coefficients of the $z$-power series yields finally an expression similar to that of the Dickson-Hipp type representation  in the \lev case (see \cite{IP})
\begin{equation} \la{ZDH}
Z_v(x,w)=\left(\pg(w)-\frac{w}{v}\right)\sum_{k=0}^\infty w^kW_v(x+k),\quad w\in (0,\vf_v),\; v\in (0,1],\; x\in \mathbb{N}_0
\end{equation}
(it is easy to check that this expression has $z$-transform \eqref{eq:Z-generating-fnct}).

In the special case $w=1$ we set $Z_v(x):=Z_v(x,1)$, \eqref{eq:Z-generating-fnct} simplifies to
\begin{equation} \label{eq:Z-basic-transform}
\T Z_v(z):=\sum_{x=0}^\infty z^xZ_v(x)=\fr{  \T p(z)-z}{ (\T p(z)-\frac{z}{v})(1-z)},\quad z\in (0,\vf_v),\; v\in (0,1],\;
\end{equation}
and we have the representation
\begin{equation}\label{eq:Z-basic}
Z_v(x)=1+\left(\frac{1}{v}-1\right)\sum_{y=0}^{x-1}W_v(y),\quad v\in (0,1],\; x\in \mathbb{N}_0.
\end{equation}

\beR
Using \eqref{eq:ruin-transform} in the form
\begin{equation*}  \la{PKdr} \T{\Rui}_v(z)=\fr{1 }{\frac{z}{v}- {\pg(z)}}\le(\fr{z-\pg(z)}{(1-v)(1-z)}+\fr{\vf_v}{v(1-\vf_v)}\ri), \quad v, z \in (0,1),\; z\ne\vf_v,\end{equation*}
it follows from \eqref{eq:Z-basic-transform} and \eqref{Wdisc} that $${\Rui}_v(x):=\sum_{n=0}^\I  v^n \Rui(n;x)=\frac{1}{1-v}Z_v(x)-\frac{\vf_v}{v(1-\vf_v)}W_v(x),$$
i.e.
\begin{equation}\label{eq:levy-identity}
E_x[v^{\tau_{-1}^-};\tau_{-1}^-<\infty]=
Z_v(x)-\frac{\vf_v(1-v)}{v(1-\vf_v)}W_v(x)=Z_v(x)-\a_v W_v(x),\quad x\in\mathbb{N}_0,\; v\in (0,1),
\end{equation}
where we have set $\a_v:=\a_v(1)$ (recall that we have chosen $\a_v(1)$ so that $Z_v(0)=1=E[v^{\tau_{-1}^-};\tau_{-1}^-<\infty]+\a_v(1)W_v(0)$). Passing to the limit $v\uparrow 1$, we find that $P_x(\tau_{-1}^-<\infty)=1-W(x)(1-\pg'(1-)\land 1)$.
\eeR

\beR
It is seen from \eqref{eq:Z-basic}, Remark~\ref{remark:Ws-connection} and \cite[Def.~4.9]{vidmar2013fluctuation} that one has the identification $Z^{(q)}(mh)=Z_{\frac{\gamma}{\gamma+q}}(m)$ for $q\in [0,\infty)$, $m\in \Z$, where $Z^{(q)}$ is the $Z$ $q$-scale function of \cite{vidmar2013fluctuation}. Then \eqref{eq:Z-generating-fnct}, \eqref{eq:Z-via-deficit-at-ruin} and \eqref{eq:Z-recurrence}, with $w=1$, become \cite[Eq.~(4.9), Prop.~4.13 and Eq.~(4.11)]{vidmar2013fluctuation}, respectively; \eqref{eq:levy-identity} becomes \cite[Eq.~(4.8)]{vidmar2013fluctuation}. For an alternative form of \eqref{eq:Z-recurrence} (when $w=1$) see \cite[Eq.~(4.14)]{vidmar2013fluctuation}.
\eeR

\beP\label{proposition:martingale-Z}
For each $v\in (0,1]$, $w\in (0,1]$, the process $(v^{n\land \tau_{-1}^-}Z_v(X_{n\land \tau_{-1}^-},w))_{n\in \mathbb{N}_0}$ is a martingale.
\eeP

\begin{proof} This follows for instance by linearity, from Proposition~\ref{prop:scales:miscellaneous}, and from the definition of $Z_v(\cdot,w)$ via the Markov property and the terminal time property of $\tau_{-1}^-$.
\end{proof}

\beR
For the case $w=1$, the analogue of Proposition~\ref{proposition:martingale-Z} 
in the setting of upwards skip-free L\'evy chains are the martingales, for $q\in [0,\infty)$, $(e^{-q(t\land \tau_{-h}^-)}Z^{(q)}(Y_{t\land \tau_{-h}^-}))_{t\in [0,\infty)}$ \cite[Cor.~4.17]{vidmar2013fluctuation}. In the case of a spectrally negative L\'evy processes $U$,  $(e^{-q(t\land \tau_0^-)}Z^{(q)}(U_{t\land \tau_0^-}))_{t\in [0,\infty)}$ is a local martingale with localizing sequence $(\tau_n^+)_{n\in \mathbb{N}}$ \cite[Ex.~8.12]{Kyp}. See also \cite{APP15}: There, Gerber-Shiu functions are defined as solutions to martingale problems \cite[Def.~5.1]{APP15}, and the $Z^{(q)}(\cdot,\theta)$ function is the Gerber-Shiu function with  boundary condition $e^{x \theta}$ for $x\in (-\infty,0)$ \cite[Def.~5.8]{APP15}.
\eeR

\beR \la{r:deficit-at-ruin} Assume $E \claim_1<\infty$; let $v\in (0,1]$, $x\in \mathbb{Z}$. We can obtain the expected undershoot at ruin by differentiating \eqref{eq:Z-via-deficit-at-ruin} with respect to $w$ from the left at $1$. Putting
$Z_{1,v}(x):=-\fr{\partial Z_{v}(x,w)}{\partial w}\vert_{w=1-}$, we find that for $b\in \mathbb{N}_0$,
\begin{equation}\label{eq:deficit-at-ruin}
E_x [X(\t^-_{-1}) \; v^{ \t^-_{-1}}   ; \t^-_{-1} <\t_b^+]=Z_{1,v}(x)-\frac{W_v(x)}{W_v(b)}Z_{1,v}(b) , \quad x\leq b.
\end{equation}
The generating function transform of $Z_{1,v}$ is given by 
\beq \la{Z1tr} \T Z_{1,v}(z) :=\sum_{k=0}^\infty z^k Z_{1,v}(k)= \fr{z}{1-z} \fr 1{\pg(z)-z/v}
\left(\frac{\T p(z)-z}{1-z}-(1-\T p'(1-))\right),\quad z\in (0,\vf_v).
\eeq
Setting for  $f:\mathbb{N}_0\to \mathbb{R}$ and $y\in \mathbb{N}_0$, $\bar f(y):=\sum_{z=0}^{y-1}f(z)$ (in particular, $\bar f(0)=0$), and using $\sum_{k=0}^\infty z^k\overline{f}(k)=\frac{z}{1-z}\sum_{k=0}^\infty z^kf(k)$ for $z\in (0,1]$, we find that for $x\in \mathbb{N}_0$, this coincides with the generating function of $\mathbb{N}_0\ni x\mapsto \bar Z_{v}(x)-(1-\T p'(1-))\bar W_{v}(x)$, i.e. 
\begin{equation}
Z_{1,v}(x)=\bar Z_{v}(x)-(1-\T p'(1-))\bar W_{v}(x), \quad x\in \mathbb{N}_0,
\end{equation}
Note also that  when $x <0$, $Z_{1,v}(x)=x$. $Z_{1,v}$ will play a central role in the modified de Finetti problem -- see Subection \ref{sec:deficit-at-ruin-refl}, and in its doubly reflected variant presented in Subection \ref{sec:combination-dividends-injections}. 
\eeR

\section{Applications to the study of a company's capital surplus process}\label{section:application}
In this section we investigate various forms of the (combined) capital injections-dividend payouts-penalty at ruin problem. One typically has in mind an insurance company, but this need not be the case.

\subsection{The moment generating function of cumulative capital injections} \la{sec:ci}
For the simplest case, we begin by considering a company, whose surplus capital process $\T X=(\T X_k)_{k\in \mathbb{N}_0}$ obeys the following dynamics: for $k\in \mathbb{N}_0$, given that at the end of period $k$, its capital is $\T X_k$, then in period $k+1$ the company receives (the premium) $1$, pays out the (claim) amount $\claim_{k+1}$, and, should its net capital at this point be strictly negative, receives a capital injection that just brings its capital back to zero at the end of the $(k+1)$-th period, i.e. $\T X_{k+1}=(\T X_k+1-L_{k+1})\lor 0$. If the initial capital $x\in \mathbb{Z}$ of the company is strictly negative, the company receives immediately the capital injection $-x$, so that its capital at the end of the zeroth period is nonnegative, i.e. $\T X_0=(-x)\lor 0$. One says that the surplus process has the dynamics of $X$ reflected at $0$.

Let then $R_*(n):=(-\inf_{m\leq n} X_m)\lor 0$, $n\in \mathbb{N}_0$, denote the cumulative capital injections for the process $X$ reflected at $0$, and let, for $b\in \mathbb{N}_0$, $\T \t_b^+$ denote the first entrance time into $[b,\infty)$ by the reflected process. It was discovered by \cite{IP} that their joint \mgf is very simply expressible in terms of the second scale function of two parameters. In our context,  their formula becomes
\beP
For $b\in \mathbb{N}_0$,
\beq \la{capinj}  B_v^b(x,w):=E_x[v^{\T \t_b^+} w^{R_*(\T \t_b^+)};\T \t_b^+ <\infty]=
\begin{cases}
\fr{Z_v(x,w)} {Z_v(b,w)} & x\leq b\\
1& x>b
\end{cases},\quad \{v,w\}\subset (0,1]. \eeq
\eeP
\begin{proof}
The case $x>b$ is trivial; assume $x\leq b$. Then this formula is ``equivalent" to  \eqref{eq:Z-via-deficit-at-ruin}, since by the strong Markov property of $X$,\footnotesize
$$E_x[v^{\T \t_b^+} w^{R_*(\T \t_b^+)};\T \t_b^+ <\infty]=
\E_x \left[v^{\t_{-1}^-} w^{-X(\t_{-1}^-)} ;  \t_{-1}^- < \t_b^+ \right]  \E_0[v^{\T \t_b^+} w^{R_*(\T \t_b^+)};\T \t_b^+ <\infty]+ \E_x\left[v^{\tau_b^+};\tau_b^+<\tau_{-1}^-\right],$$\normalsize
i.e.
\begin{equation} \la{proof} B_v^b(x,w)=
\E_x \left[v^{\t_{-1}^-} w^{-X(\t_{-1}^-)} ;  \t_{-1}^- < \t_b^+ \right]  B_v^b(0,w)+ W_v(x)W_v(b)^{-1}.\end{equation}
Thus, if $ B_v^b(x,w)$ is known from \eqref{capinj}, one gets an equation for the deficit at ruin quantities
\begin{eqnarray*}&& Z_v(x,\tw) Z_v(b,\tw)^{-1}= W_v(x)W_v(b)^{-1} + \E_x \left[v^{\t_{-1}^-} w^{-X(\t_{-1}^-)} ; \t_{-1}^- < \t_b^+ \right] Z_v(b,\tw)^{-1},  \end{eqnarray*}
with  solution \eqref{eq:Z-via-deficit-at-ruin}.
 And if the solution to the deficit at ruin problem is known as \eqref{eq:Z-via-deficit-at-ruin}, one may use \eqref{proof} to obtain, first with $x=0$, $B_v^b(0,w)=Z_v(b,\tw)^{-1}$, and then \eqref{capinj}.
\end{proof}

\subsection{The  \deF dividends optimization problem}
\label{sec:deF}
Now the company pays dividends, but does not receive capital injections.
Letting for $k\in \mathbb{N}_0$, $r(k)$ denote the dividend amount (necessarily $\mathbb{N}_0$-valued)
paid out at the end of period $k$, we have the following dynamics for the end-of-period surplus process $\tilde{X}$: for $k\in \mathbb{N}$, in period $k$, the company receives $1$, pays out $\claim_k$ and then, assuming ruin has not yet occurred, the amount $r(k)$, yielding
$\tilde{X}_k=\tilde{X}_{k-1}+1-\claim_k-r(k)$. Once ruin has occurred, the process is stopped, and no dividends are paid out thereafter. At end of period zero, if the initial capital $x\in \mathbb{Z}$ is strictly positive, the dividend amount $r(0)$ is paid out, so that $\tilde{X}_0=x-r(0)$. We insist $r(k) \leq \tilde{X}_{k-1}+1-\claim_k$ for $k\in \mathbb{N}$ and $r(0)\leq x$ (i.e. dividend payouts cannot lead to ruin). The dividend policy process $(r(k))_{k\in \mathbb{N}_0}$ must be adapted to the natural filtration of $(\claim_k)_{k\in \mathbb{N}}$.

The classic \deF problem then consists in computing the optimal
discounted dividends until ruin under all  dividend policies satisfying the above constraints -- see de Finetti \cite{DeF},
 Miller and Modigliani   \cite{miller1961dividend}
 (in a deterministic setup), Miyasawa \cite{miyasawa1961economic}
 and Gerber \cite{gerber1972games}. Here we agree that in the optimization objective, $r(k)$ is discounted (multiplied) by $v^k$, where $v\in (0,1]$ is the discount factor. To exclude some degeneracy, we assume throughout this subsection that $(p_0+p_1)\land v<1$.
 \beD \la{d:bp}
For $b\in \mathbb{N}_0$, a  dividend  policy $\pi_b$ with  barrier $b$ consists in taking $r(0)=(x-b)^+$ and $r(k)=(\tilde{X}_{k-1}+1-\claim_k-b)^+$ for $k\in \mathbb{N}$, up to ruin, i.e. (since we are in the upwards skip-free case) in reducing
the reserves  each time they reach $b+1$ (except possibly at time zero, when $x-b$ may be strictly larger than $1$). We will write the expectation operator $\Eb_x$ and the probability $\P^{b]}_x$ to indicate this policy and the initial capital $x$. One says that under $\Eb_x$, $\tilde{X}$ follows the dynamics of the process $X$ reflected at $b$. The sets
$$\mC^b:=[0,b]\text{ and } \mD^b:=(b,\I)$$
are called the continuation and dividend taking set, respectively; $\T\tau^+_b:=\inf\{k\in \mathbb{N}_0:\T X_k\geq b\}$.

The ruin time, i.e. the first time the surplus process becomes strictly negative, will be denoted by $\Tt$. Note that $r(k)=0$ for $k\geq \Tt$. We also set, for $k\in \mathbb{N}_0$, $R(k):=\sum_{i=0}^{k}r(i)$, the cumulative dividends paid out up to (including) period $k$, and interpret $R(k)=0$ for $k<0$.
\eeD

\beP
The value function
under a barrier dividend distribution policy
$\pi_b$ with barrier $b\in \mathbb{N}_0$
 is given by:
\beq\label{e:Dobj}
  V^b_{D}(x):=\Eb_x \sum_{i=0}^{\Tt-1} v^i r(i)=\bc
  \fr{W_{v}(x)}  {\D W_{v}(b)} & x \leq b\\
  x-b + V^b_{D}(b)=x-b+\fr{W_v(b)}{\D W_v(b)} & x > b\ec, \eeq
where for
$f:\mathbb{N}_0\to \mathbb{R}$, $k\in \mathbb{N}_0$,
$\D f(k):=f(k+1)-f(k)$ gives the forward difference operator.
\eeP
 \beR It is clear from \eqref{Wdef} that under the stipulation $(p_0+p_1)\land v<1$, $W_v$ is strictly increasing. \eeR

\begin{proof} The case $x>b$ is trivial; assume $x\leq b$. Then \eqref{e:Dobj} is ``equivalent" to \eqref{Wdef}, since using the strong Markov property of $X$, one has clearly the relation
\begin{equation}\label{proof:divs}
V^b_D(x)=\E_x[v^{\tau_{b+1}^+};\tau_{b+1}^+<\tau_{-1}^-](1+V^b_D(b)).
\end{equation}
Thus if \eqref{e:Dobj} is known, one obtains \eqref{Wdef}, and vice versa, if \eqref{Wdef} is known, then one obtains by setting $x=b$ in \eqref{proof:divs}, first $V^b_D(b)=\frac{W_v(b)}{\Delta W_v(b)}$ and then by substituting back, \eqref{e:Dobj}.
\end{proof}

  \beR \label{remark:optimization-dividends}The ``factorization result" $\fr{W_{v}(x)}  {\D W_{v}(b)}$  of \eqref{e:Dobj} has been known for a long time \cite[Eq.~(19)]{morrill1966one}  \cite[Sec. 5, Eq.~(3.1)]{gerber2010elementary}, and in the simplest case when $\D W_v$  is ``unimodal with minimum at $b^*$'', i.e. when $\D W_v$ is nondecreasing after $b^*$ and nonincreasing before $b^*$, it yields in fact the optimal value function over all dividend distribution policies. The optimal ``barrier"  policy of taking dividends
 in $\mD^{b^*}=(b^*,\I)$ and continuing in $\mC^{b^*}=[0,b^*]$ can then be viewed as a
 transformation of the \sf into the value function
$V^b_D$, which must be concave, by  ``linearization" of the
 convex piece of $W_v$.\footnote{By a convex (concave) function
 $f:\mathbb{N}_0\to \mathbb{R}$ we mean a function whose
 forward difference $\D f$ is nondecreasing (nonincreasing).}
 \eeR
 When $\D W_v$  is not unimodal, the optimal policy may be ``multi-band", and requires a complicated recursive construction \cite{morrill1966one,Schmidli,APP15}. We will recall this concept briefly in Definition \ref{def:multiband}, but the main concern of our applications is optimization among barrier policies, by which we mean optimizing the limit $b$ of the continuation interval $[0,b]$, in the sense of finding
$$ V_D(x):=\sup_{b\in \mathbb{N}_0}V^b_D(x).$$ With the objective given by \eqref{e:Dobj}, this is related to maximizing the ``barrier influence function'' $1/\D W_v$, i.e. minimizing $\D W_{v}$ (as is customary\footnote{And we will follow an analogous convention with respect to the optimization problems of Subsections~\ref{sec:deficit-at-ruin-refl}
 and~\ref{sec:combination-dividends-injections} to follow.}, we will say $b\in \mathbb{N}_0$ is optimal for $V_D(x)$, if $V_D(x)=V^b_D(x)$):

\beL  \label{lemma:optimization-dividends}(I) If $q:=\inf_{b\in \mathbb{N}_0}\Delta W_v(b)$ is attained, letting $b^*$ be any minimizer of $\Delta W_v$, it follows that  $b^*$ is optimal for $V_D(x)$, whenever $x \leq b^*$. (II)  If  the infimum defining $q$ is not attained, then the supremum defining $V_D(x)$ is not attained either. (III)  If for some $b\in \mathbb{N}_0$, the function $\Delta W_v$ is nondecreasing after $b^*$, i.e. satisfies $\Delta W_v(b) \leq \Delta W_v(b')$ for all $b' > b \geq b^*$,  and nonincreasing before $b^*$, i.e. satisfies $\Delta W_v(b) \geq \Delta W_v(b')$ for all $b < b '\leq   b^*$, and if furthermore $b^*<x$, then $b^*$ is optimal for $V_D(x)$.
\eeL

\beR
This dovetails nicely with Remark~\ref{remark:optimization-dividends}: when $\Delta W_v$ is unimodal with minimum at $b^*$, then $b^*$ is optimal for $V_D(x)$, \emph{whether or not} $x\leq b^*$.
\eeR

\begin{proof}
(I) To see this, note that for $ x \leq b$, $V^b_D(x)=\frac{W_v(x)}{\Delta W_v(b)}\leq \frac{W_v(x)}{\Delta W_v(b^*)}=V^{b^*}_D(x)$. And for $b<x$, $V^b_D(x)=x-b+\frac{W_v(b)}{\Delta W_v(b)}\leq x-b+\frac{W_v(b)}{\Delta W_v(b^*)}\leq \frac{W_v(x)}{\Delta W_v(b^*)}=V^{b^*}_D(x)$, where the final inequality follows from (telescopic sum) $W_v(x)-W_v(b)=\sum_{k=b}^{x-1}\Delta W_v(k)\geq \sum_{k=b}^{x-1}\Delta W_v(b^*)=(x-b)\Delta W_v(b^*)$. (II) Indeed, there exists a sequence $(b_n)_{n\in \mathbb{N}}$ in $\mathbb{N}$, with $\Delta W_v(b_n)$ satisfying $\Delta W_v(b')>\Delta W_v(b_n)$ for all $b'< b_n$, $n \in \mathbb{N}$. Let now $b\in \mathbb{N}_0$. There is an $n\in \mathbb{N}$ such that $b_n\geq x\lor b$. Then if $b\geq x$, clearly $V^{b}_D(x)=\frac{W_v(x)}{\Delta W_v(b)}\leq \frac{W_v(x)}{\Delta W_v(b_n)}=V^{b_n}_D(x)$. And if $b<x$, then $V^b_D(x)=x-b+\frac{W_v(b)}{\Delta W_v(b)}\leq x-b+\frac{W_v(b)}{\Delta W_v(b_n)}\leq \frac{W_v(x)}{\Delta W_v(b_n)}=V^{b_n}_D(x)$, where the last inequality follows from $W_v(x)-W_v(b)=\sum_{k=b}^{x-1}\Delta W_v(k)\geq \sum_{k=b}^{x-1}\Delta W_v(b_n)=(x-b)\Delta W_v(b_n)$.  In other words, as $n\uparrow \infty$, $V^{b_n}_D(x)\uparrow\uparrow \sup_{b\in \mathbb{N}_0}V^b_D(x)=V_D(x)$, which however is not attained. We also see that $q>0$, since  $V_D^b(x)$ is bounded by $(x-b)^++\sum_{k=1}^\infty v^k$, as $b$ ranges over $\mathbb{N}_0$. (III) Since for $y\in \mathbb{N}_0$, $\left(\frac{W_v(y+1)}{\Delta W_v(y+1)}-(y+1)\right)-\left(\frac{W_v(y)}{\Delta W_v(y)}-y\right)=W_v(y)\left((\Delta W_v(y+1))^{-1}-(\Delta W_v(y))^{-1}\right)$, it follows from the assumption, that the map $\mathbb{N}_0\ni b\mapsto \frac{W_v(b)}{\Delta W_v(b)}+(x-b)$ has a maximum at $b^*$. Thus if $b\leq x$, then it follows at once that $V^b_D(x)\leq V^{b^*}_D(x)$. And if $b> x$, then $ V^{b^*}_D(x)\geq  V^{x}_D(x)=\frac{W_v(x)}{\Delta W_v(x)}\geq \frac{W_v(x)}{\Delta W_v(b)}=V^b_D(x)$.
\end{proof}

\beR
For $x\leq b$, by the skip-free property, \footnotesize
\beq \la{o=g} V_D^b(x)=\Eb_x\left[\sum_{n=1}^\infty v^n \mathbbm{1}(n < \Tt,\T X_{n-1}=b, \claim_n =0) \right]=\Eb_x\left[\sum_{n=1}^{\Tt \wedge \kil_v - 1} \mathbbm{1}( \T X_{n-1}=b, \claim_n =0)\right]=\Eb_x R(\Tt \wedge \kil_v - 1).\eeq\normalsize
where $\kil_v$ is an independent random variable with distribution $\geom_\mathbb{N}(1-v)$.\footnote{For $r\in (0,1]$, we denote by $\geom_\mathbb{N}(r)$, resp. $\geom_{\mathbb{N}_0}(r)$, the geometric law on $\mathbb{N}$, resp. $\mathbb{N}_0$, with success parameter $r$, i.e. having p.m.f. $\mathbb{N}\ni k\mapsto r(1-r)^{k-1}$, resp.  $\mathbb{N}_0\ni k\mapsto r(1-r)^{k}$. The degenerate cases $\geom_\mathbb{N}(0)$ and $\geom_{\mathbb{N}_0}(0)$ are both interpreted as $\delta_\infty$, the Dirac mass at $\infty$.}
\eeR
\beXa For $b=0$, plugging  $W_v(0)=p_0^{-1}$ and $W_v(1)=p_0^{-2}(v^{-1}-p_1)$ into \eqref{e:Dobj}, yields
\beq \la{gc} V_D^0(0)= \fr{p_0 v}{1 -p_1 v -p_0 v}.\eeq

For $p_1=0$, this reduces to  (note that we start with initial capital zero, hence pay no dividends at time zero, and  that   dividends of $1$ are taken all the times strictly prior to ruin) \beq V_D^0(0)=
\fr{p_0 v}{1- p_0v}=\Enull_0\left[\sum_{n=1}^\infty v^n \mathbbm{1}_{\{n < \Tt\}}\right]=\Enull_0[\Tt \wedge \kil_v - 1]=\Enull_0 R(\Tt \wedge \kil_v - 1),
\eeq where $\Tt\sim\geom_\mathbb{N}(1-p_0)$ and  $\kil_v\sim\geom_\mathbb{N}(1-v)$ and hence $R(\Tt \wedge \kil_v-1)=\Tt \wedge \kil_v-1\sim \geom_{\mathbb{N}_0}(1-p_0v)$.

When $p_1$ is not necessarily equal to $0$, one may still decompose $X$ into the process which records $X$ only when it changes its value  --- it does so each time independently according to the law of $\claim_1$ conditioned on $\{\claim_1\ne 1\}$ ---  and into the independent amounts of time that elapse in-between these changes, them being i.i.-$\geom_\mathbb{N}(1-p_1)$-d.\footnote{This  is analogous to the decomposition of a continuous time Markov chain into its jump chain and its sojourn times.} From the perspective of the surplus process, this means that it may be seen as evolving (up to ruin) according to the following probabilistic prescription: for $k\in \mathbb{N}_0$, if $0$ at end of period $k$ (i.e. ruin has not yet occurred), then for $L$ subsequent periods, where $L\sim\geom_{\mathbb{N}_0}(1-p_1)$, the claims are equal to $1$, just off-setting the premia, and then during period $k+L$, independently, the surplus process goes up by $1$ with probability $p_0/(1-p_1)$ or down by $l$ with probability $p_{l+1}/(1-p_1)$, $l\in \mathbb{N}$ -- if the former, a dividend of one is taken; if the latter, ruin occurs. It follows that in this case the total discounted dividends are equal to $$V^0_D(0)=\E \left[\sum_{i=1}^{\tau-1}v^{\sum_{j=1}^iQ_j}\right]=\E \sum_{i=1}^{\tau-1}\left(\frac{v(1-p_1)}{1-vp_1}\right)^i=
\frac{\frac{p_0}{1-p_1}\frac{v(1-p_1)}{1-vp_1}}{1-\frac{p_0v}{1-vp_1}}
=\frac{p_0v}{1-p_1 v-p_0v}, $$ where $\tau\sim\geom_\mathbb{N}(1-\frac{p_0}{1-p_1})$ and $Q_j\sim \geom_\mathbb{N}(1-p_1)$, $j\in \mathbb{N}$, are independent, confirming again \eqref{gc}. In other words, it is the same as the case $p_1=0$, except that one has conditioned the claims not to be equal to $1$, $p_0\rightsquigarrow \frac{p_0}{1-p_1}$, and changed the discount factor, $v\rightsquigarrow \frac{v(1-p_1)}{1-vp_1}$, reflecting the $\geom_\mathbb{N}(1-p_1)$ distributed ``holding periods'' during which $X$ does not move. Thus, for all intents and purposes, the case $p_1\ne 0$ is reduced to the case $p_1=0$. For instance, under $P^{0]}_0$, the law of the cumulative paid-out dividends, i.e. of $R(\Tt-1)$, is $\geom_{\mathbb{N}_0}(1-\frac{p_0}{1-p_1})$, and hence
\beq \la{mgp1} R(\Tt\land \kil_v-1)\sim\geom_{\mathbb{N}_0}\left(1-\frac{p_0v}{1-p_1v}\right)\eeq  (replacing $p_0$ and $p_1$ by $p_0v$ and $p_1v$, respectively, has the same effect as independent geometric killing with probability $1-v$ (the mass $(1-v)(p_0+p_1)$ may, for instance, be added to $p_2$, it matters not)). See Proposition~\ref{proposition:geometrics} below for a generalization.

Finally, expanding \eqref{gc} in $v$-series, reveals that the probability that dividends are paid in the $n$-th step is
$$\P^{0]}_0[r(n)= 1]=\P^{0]}_0[\Tt >n,\claim_n =0]=(p_0+ p_1)^{n-1}p_0, \quad n\in \mathbb{N},$$ which also has a clear interpretation: $(n-1)$-times ruin must not occur, i.e. the claim is zero or one, and then the $n$-th claim must be zero.  Incidentally, the above is the survival function of a  modified \geo r.v. $\T T$ with
\begin{equation*} \la{mg} P[\T T=1]=1-p_0,\; P[\T T=k]=p_0 (1-p_0-p_1)(p_0+p_1)^{k-2},\quad k \in \mathbb{N}_{\geq 2}.\end{equation*}
\eeXa

The next result gives another probabilistic interpretation to the
objective $V_D^b(b)=\frac{W_v(b)}{\Delta W_v(b)}=\frac{W_v(b+1)}{\Delta W_v(b)}-1=\left(\frac{\Delta W_v(b)}{W_v(b+1)}\right)^{-1}-1$, which is the mean of $\geom_{\mathbb{N}_0}\left(\frac{\Delta W_v(b)}{W_v(b+1)}\right)$. Note that much more is known in  the case of \sn \lev processes, where $\le(V_D^b(b)\ri)^{-1}=\frac{(W^{(q)})'(b)}{W^{(q)}(b)}$, coincides with the  rate of ``excursions" larger than $b$ of the  Poisson process of heights of downward excursions from a
running maximum, in the presence of exponential killing at rate $q$ -- see \cite[Sec.~VII.8]{Ber} for $q=0$ and \cite{Doney} for $q>0$.

\beP\label{proposition:geometrics} Let $b\in \mathbb{N}_0$. Under a barrier policy $\pi_b$, starting from $x=b$, the killed cumulative dividends until ruin, $R(\Tt \wedge \kil_v - 1)$, have the law $\geom_{\mathbb{N}_0}\left(\frac{\Delta W_v(b)}{W_v(b+1)}\right)$ (recall $\kil_v\sim\geom_\mathbb{N}(1-v)$, independent of $\T X$). In particular, \beq \la{mgfR} \Eb_b z^{R(\Tt \wedge \kil_v -1)}=\fr{1-
\frac{ W_v(b)}{W_v(b+1)}}{1-z \frac{ W_v(b)}{W_v(b+1)}},\quad z\in (0,1].\eeq
\eeP

\begin{proof}
First one assumes $p_0+p_1<1$ and $v=1$. We have the representation of $R(\Tt \wedge \kil_v - 1)$ as the sum $\sum_{i=1}^N \T R_i$, where $\T R_i$ are i.i.d. with the law given in \eqref{mgp1}, and $N\sim \geom_\mathbb{N}(1-\a(b))$ is an independent geometric r.v. with $\a(b)$ yet to be determined. Indeed, the successive $\T R_i$ come from the dividends collected during the periods of time that the surplus process either stays at the level $b$, or else increases to $b+1$, only to be taken down to $b$ by a paid-out dividend. These amounts have the same law as does the amount of dividends collected until ruin when starting from $0$ under $\pi_0$. On the other hand, $\a(b)$ is the probability that the surplus process, once it has jumped to a level strictly below $b$, then goes on to reach the level $b$ before ruin occurs, i.e. (the quotients $\frac{p_k}{1-p_0-p_1}$ come from conditioning to jump strictly below $b$ from $b$) $\a(b)=\sum_{k=2}^{b+1}\frac{p_k}{1-p_0-p_1}\P_{b-k+1}[\tau_b^+<\tau_{-1}^-]=\sum_{k=2}^{b+1}\frac{p_k}{1-p_0-p_1}\frac{W(b-k+1)}{W(b)}$, which equals, using \eqref{har}, $\frac{W(b)-p_1W(b)-p_0W(b+1)}{(1-p_0-p_1)W(b)}=1-\frac{p_0}{1-p_0-p_1}\frac{\Delta W(b)}{W(b)}$. The conclusion of the proposition then follows e.g. by computing the probability generating function of the ``geometric sum of geometrics" $\sum_{i=1}^N \T R_i$ and recognizing the geometric random variable and its parameter. The general case for $p_0+p_1<1$  is got by replacing $p_0,\ldots,p_{b+1}$ by $p_0v,\ldots,p_{b+1}v$ (and for instance adding the mass $(1-v)(p_0+\cdots+p_{b+1})$ to $p_{b+2}$, it matters not), using Remark~\ref{remark:killing}. When $p_0+p_1=1$, then the result clearly still holds true (one gets, using \eqref{har}, the law of \eqref{mgp1}, i.e. $\geom_{\mathbb{N}_0}(\frac{1-v}{1-p_1v})$, as one should).
\end{proof}

The following proposition gives a dividends-deficit at ruin type law for the compound binomial risk processes reflected at $b$, in the style of \cite[Sec.~4]{gerber2010elementary}. See \cite[Thm.~6]{IP}, \cite[Lem.~6]{AGV} for the \lev analog.

\beP \la{c:DP}  The   joint generating function of the ruin time, deficit at ruin and of the cumulative dividends for a compound binomial risk process reflected at $b\in \mathbb{N}_0$  is given by, with $\{v,z,w\}\subset (0,1]$,
\begin{equation} \la{e:DP}
 DP_v^{b]}(x,w,z):=\Eb_x\left[v^{  \Tt}  w^{-\T X(\Tt)}z^{R(\Tt)};\Tt<\infty \right]=
\begin{cases}
Z_v(x,w) - \frac{Z_v(b+1,w)-zZ_v(b,w)}{W_v(b+1)-zW_v(b)}W_v(x) & x\leq b\\
z^{x-b}DP_v^{b]}(b,w,z)& x>b
\end{cases}.
\end{equation}
\eeP
\beR
When $p_0+p_1<1$, by setting $v=z=w=1$, one obtains (as one should) $\P^{b]}_x(\Tt<\infty)=1$ for all $x\in \mathbb{Z}$. When $p_0+p_1=1$, we have of course $\Tt=\infty$, $\P_x^{b]}$-a.s. for all $x\in\mathbb{N}_0$ (and $\Tt=0$, $\P_x^{b]}$-a.s. for all $x\in -\mathbb{N}$).
\eeR
\begin{proof}
The case $x>b$ is trivial; let $x\leq b$.  Using the strong Markov property for $X$ at
the exit time from the interval $[0,b)$ yields that $g(x):=DP_v^{b]}(x,w,z)$ satisfies:
\begin{eqnarray*}
 &&g(x)=Z_v(x,w) -     \fr{W_v(x)   }{ W_v(b)} Z_v(b,w)+ \fr{ W_v(x)   }{ W_v(b)} g(b)=Z_v(x,w) + W_v(x)\fr{ g(b)- Z_v(b,w) }{ W_v(b)} \\&& \Rightarrow \fr{g(x)-Z_v(x,w)}{W_v(x)   } =\fr{g(b)-Z_v(b,w)}{W_v(b)}=:-H_v^b(z,w).
\end{eqnarray*}
Now by conditioning on the first jump
$$g(b)= p_0vzg(b)+ p_1 v g(b)+ \sum_{k=1}^{b}p_{k+1}vg(b-k)+v\sum_{k=b+1}^\infty p_{k+1}w^{k-b}.$$
Plugging in $g(x)=Z_v(x,w)-W_v(x)H_v^b(z,w)$ gives us \footnotesize
$$(1-p_0vz-p_1v)(Z_v(b,w)-W_v(b)H_v^b(z,w))=\sum_{k=1}^{b}p_{k+1}v(Z_v(b-k,w)-W_v(b-k)H_v^b(z,w))+v\sum_{k=b+1}^\infty p_{k+1}w^{k-b}.$$\normalsize
Using now \eqref{har} and \eqref{eq:Z-recurrence} reduces this to
$$(1-p_0vz-p_1v)(Z_v(b,w)-W_v(b)H_v^b(z,w))=$$
$$-H_v^b(z,w)(W_v(b)-vp_1W_v(b)-vp_0W_v(b+1))+Z_v(b,w)-vp_1Z_v(b,w)-vp_0Z_v(b+1,w),$$  i.e. $p_0v(Z_v(b+1,w)-zZ_v(b,w))=p_0vH_v^b(z,w)(W_v(b+1)-zW_v(b))$.
 \end{proof}
\beR
As a check, setting $v=w=1$ and $x=b$ in \eqref{e:DP} recovers \eqref{mgfR} in the case $v=1$.
\eeR

Taking $z=1$ in \eqref{e:DP} yields
\beC \la{c:rt}
For $\{v,w\}\subset (0,1]$, the   joint generating function of the (reflected) ruin time and of  the deficit at ruin for a compound binomial risk process reflected at $b\in \mathbb{N}_0$  is given by
\begin{equation} \la{e:DP1}
 \Rui_v^{b]}(x,w):=\Eb_x\left[v^{  \Tt}  w^{-\T X(\Tt)};\Tt<\infty\right]=
Z_v(x,w) - \frac{\D Z_v(b,w)}{\D W_v(b)}W_v(x),  \quad x\leq b.\null\hfill\qed
\end{equation}
\eeC

\beR This result is similar to identity \eqref{eq:Z-via-deficit-at-ruin} for the joint generating function of the  ruin time and of  the deficit at ruin, with
absorbtion at $b$; this is to be expected, since we only replaced
the boundary condition $\Rui_v^{b}(b,w):=E_b [v^{ \t^-_{-1}} w^ {-X(\t^-_{-1})} ; \t^-_{-1} <\t_b^+]=0$ by $\D \Rui_v^{b]}(b,w)=0$. \eeR

We recall finally some further background information for the general  \deF dividends optimization problem with no penalty for the deficit at ruin, when $\D W_v$ is not unimodal. This is useful for the numerics Section~\ref{sec:examples}, to understand the examples where the optimal dividends policy is  ``multi-band".
 \beD\label{def:multiband} A multi-band dividends policy is specified by
  a partition of $\N$ into  continuation intervals $\mC_1=[0,b_1]$, $\mC_2=[a_2,b_2]$, \ldots,
 and dividend taking intervals $\mD_1=(b_1,a_2)$, $\mD_2=(b_2,a_3)$,
 \ldots, {intertwined} as follows:
 $\mC_1 < \mD_1 < \mC_2 < \mD_2<\ldots$.
 When the {capital} position is  in $\mD_i$,
 dividends are taken bringing the process down to the upper boundary
 $b_i$ of  $\mC_i$.
 \eeD

 When there is only one such pair $\mC_1=[0,b_1], \mD_1=(b_1,\I)$, this is the  barrier policy
$\pi_{b_1}$ of Definition~\ref{d:bp}. Subsequent $\mC_i$ and $\mD_i$, $ i \in \mathbb{N}_{\geq 2}$, appear in the optimal policy when $\D W_v$ is not unimodal and its global minimum $b_1$ is followed by other local minima. Intuitively, the existence of local minima succeeding the global one offers incitement to postpone bringing the process to $b_1$ (and thus the eventual ruin below $-1$) -- see \cite{morrill1966one} for more details.\footnote{The barriers  $b_i$, $i\in \mathbb{N}_{\geq 2}$,  may arise then, by ``shifting optimally" these local minima. See \cite{APP15} for a  recursive algorithm achieving this, which is based on the idea that the process starting in $\mC_i$  will never visit states above
$b_i +1$. Since the process at $x$ only needs to see the bands below $x$,
 $b_1$ may be computed as if only barrier policies were allowed,
i.e. taken at the global maximum of the barrier influence function. For $[a_2,b_2]$, however,
we need to take into account that the process may jump down either to ruin, or into
$\mC_1 \cup \mD_1$. Now the latter case can be viewed as termination with final payoff given by the value function $V_D^{b_1}$ over barrier policies, and this allows computing a value function $V_D^{b_1,a_2,b_2}$, and so on.}

\beR
The first
multi-band example
 is  \cite[Ex. 2]{morrill1966one}, and in the \lev case \cite{AM05}; also, the absence of local minima after the global one is known to be sufficient
  for the optimality of single barrier policies, and
   sufficient conditions in terms of the \lev measure
 have  been provided  in
   \cite[Thm. 2]{Loef}.
   However, until today,  no necessary and sufficient condition in
   terms of $W_v$ has been provided.
\eeR

\subsection{Deficit at ruin with reflection at an upper boundary and
the modified \deF problem}\label{sec:deficit-at-ruin-refl}

This problem is masterly dealt with in \cite{gerber2010elementary}.
It may be useful however to provide an alternative treatment via the
$\F,W,Z$ paradigm, as in the parallel
\lev papers \cite{loeffen2009optimal,LR,APP15,AGV}.

Specifically, we assume $\E \claim_1<\infty$ in addition to $v\land (p_0+p_1)<1$, and consider the \deF problem with dividends and no capital injections of Subsection~\ref{sec:deF}, modified by the addition of an extra linear penalty/bailout cost $k y$ upon ruin ($y$ being the (positive) deficit at ruin; $k\in (0,\infty)$). Under barrier strategies, this requires the computation of $Z_{1,v}$ \eqref{eq:deficit-at-ruin} (and of $Z_v(\cdot,w)$ \eqref{eq:Z-generating-fnct} under exponential ``risk-sensitive" bailout costs \cite{bauerle2015risk}). In precise terms, we have that, under a barrier strategy $\pi_b$, $b\in \mathbb{N}_0$, the additional expected (positive) final bailout is $k V^b_{B}(x)$, where

\beP For $b\in \mathbb{N}_0$,
    \beq \la{refldef} V^b_{B}(x):=
     \Eb_x [  v^{\Tt}  (-\T X(\Tt));\Tt<\infty]=\bc

      W_{v}(x)\fr{\D  Z_{1,v}(b)}  {\D W_{v}(b)}-Z_{1,v}(x) &x \leq b
     \\ V^b_{B}(b)& x>b \ec.
  \eeq
\eeP
\begin{proof}
In the nontrivial case, when $x\leq b$,
 using the strong Markov property for $X$ at
the exit time from the interval $[0,b)$ yields:
\footnotesize
\begin{equation}
V^b_{B}(x)=\E_x [  v^{\t^+_{b}}; \t_{b}^+ < \td]V^b_B(b)+
\E_x [  v^{\td}  (-X(\td));\td<\t_{b}^+]=\fr{W_{v}(x)}  {W_{v}(b)} V^b_{B}(b) -\le(Z_{1,v}(x) -
\fr{W_{v}(x)}  {W_{v}(b)}  Z_{1,v}(b) \ri),  \la{3F}\end{equation}\normalsize
where the second term was computed in Remark \ref{r:deficit-at-ruin}.
Making $x=0$ yields $V^b_{B}(0)=\fr{W_{v}(0)}  {W_{v}(b)} V^b_{B}(b) +\fr{W_{v}(0)}  {W_{v}(b)}  Z_{1,v}(b) $,
and substituting it back in \eqref{3F} gives us
\beq\label{eq:aux} V^b_{B}(x)+Z_{1,v}(x)= {W_{v}(x)}
\fr{V^b_{B}(0)}{W_{v}(0)}.\eeq This formula coincides with \eqref{refldef}, up to showing
that $\fr{V^b_{B}(0)}{W_{v}(0)}=
\fr{\D  Z_{1,v}(b)}  {\D W _{v}(b)}$. To see this, note that using the strong Markov property for $X$ at the exit time from the interval $[0,b]$ yields $V^b_{B}(0)=\E_0 [  v^{\t^+_{b+1}}; \t_{b+1}^+ < \td]V^b_B(b)+
\E_0 [  v^{\td}  (-X(\td));\td<\t_{b+1}^+]=\fr{W_{v}(0)}  {W_{v}(b+1)} V^b_{B}(b) +
\fr{W_{v}(0)}  {W_{v}(b+1)}  Z_{1,v}(b+1)$. Plugging into this \eqref{eq:aux} with $x=b$, i.e. $V^b_{B}(b)=-Z_{1,v}(b)+ {W_{v}(b)}
\fr{V^b_{B}(0)}{W_{v}(0)}$, we obtain the desired identity.
\end{proof}

It seems on the basis of numerics examples,
 that adding a bailout penalty typically makes
  the optimal  policy  single barrier. With this in mind and for simplicity,
  we restrict here to the version of the problem,
under which only
barrier dividend policies are allowed. Under this proviso,
optimizing under  barrier policies
the combined  objective  $$V(x):=\sup_{b\in \mathbb{N}_0} V^b(x), \quad V^b(x):=V^b_{D}(x)- k V^b_{B}(x),$$
amounts to optimizing the relevant linear combination of
the expressions \eqref{e:Dobj} and  \eqref{refldef}, viz.
$V^b(x)=(x\lor b)-b+W_v(x\land b)H(b)+kZ_{1,v}(x\land b)$,
where $H$, the ``barrier influence
 function", is given by
  \begin{equation*}
H(b):=\fr{1-k \D  Z_{1,v}(b)}  {\D W_{v}(b)}=
  \fr{1 -k\le(Z_{v}(b)-(1-\T p'(1-))W_{v}(b)\ri)}
{ \D W_{v}(b)};
\end{equation*}
see \cite[Eq.~(86)]{AGV} for the \lev case. Finding the optimum $V(x)$ is related to maximizing $H$ (cf. Lemma~\ref{lemma:optimization-dividends}):
\beL (I) If $r:=\sup_{b\in \mathbb{N}_0}H(b)$ is attained, letting $b^*$ be any maximizer of $H$, then $b^*\geq x$ implies that $b^*$ is optimal for $V(x)$. (II) If the supremum defining $r$ is not attained, then the supremum defining $V(x)$ is not attained either.
\eeL
\begin{proof}
(I) To see this, note that for $b\geq x$, $V^b(x)=W_v(x)H(b)+kZ_{1,v}(x)\leq W_v(x)H(b^*)+kZ_{1,v}(x)=V^{b^*}(x)$. And for $b<x$, $V^b(x)=x-b+W_v(b)H(b)+kZ_{1,v}(b)\leq x-b+W_v(b)H(b^*)+kZ_{1,v}(b)\leq W_v(x)H(b^*)+kZ_{1,v}(x)=V^{b^*}(x)$, where the final inequality follows from (telescopic sum) $H(b^*)(W_v(x)-W_v(b))=H(b^*)\sum_{l=b}^{x-1}\Delta W_v(l)\geq \sum_{l=b}^{x-1}(1-k\D Z_{1,v}(l))=(x-b)-k(Z_{1,v}(x)-Z_{1,v}(b))$.  (II) Indeed, there exists a sequence $(b_n)_{n\in \mathbb{N}}$ in $\mathbb{N}$, with $H(b_n)$ satisfying $H(b')<H(b_n)$ for all $b'< b_n$, $n \in \mathbb{N}$. Let now $b\in \mathbb{N}_0$. There is an $n\in \mathbb{N}$ such that $b_n\geq x\lor b$. Then if $b\geq x$, clearly $V^{b}(x)=W_v(x)H(b)+kZ_{1,v}(x)\leq W_v(x)H(b_n)+kZ_{1,v}(x)=V^{b_n}(x)$. And if $b<x$, then $V^b(x)=x-b+W_v(b)H(b)+kZ_{1,v}(b)\leq x-b+W_v(b)H(b_n)+kZ_{1,v}(b)\leq W_v(x)H(b_n)+kZ_{1,v}(x)=V^{b_n}(x)$, where the final inequality follows from $H(b_n)(W_v(x)-W_v(b))=H(b_n)\sum_{l=b}^{x-1}\Delta W_v(l)\geq \sum_{l=b}^{x-1}(1-k\D Z_{1,v}(l))=(x-b)-k(Z_{1,v}(x)-Z_{1,v}(b))$.  In other words, as $n\uparrow \infty$, $V^{b_n}(x)\uparrow\uparrow \sup_{b\in \mathbb{N}_0}V^b(x)=V(x)$, which however is not attained.
\end{proof}

\subsection{Optimizing a combination of  dividends and capital injections for a doubly reflected process}\label{sec:combination-dividends-injections}
This problem is another very good illustration of the $\F,W,Z$ paradigm and is quite hard  analytically. Indeed, the recent paper \cite{wu2011optimal} falls short of reaching an explicit  solution, which has been however available in the \lev literature \cite{APP} for a while. Since the \lev solution is a consequence of the Markov and skip-free properties, we may expect that it continues to hold in the discrete setup; and this is indeed the case.

We assume claims have a finite mean, $E \claim_1<\infty$, and linear capital injection costs $w(y)=ky$ ($y$ being the capital injection), where $k\in (1,\infty)$ is a proportionality parameter. There is also a fixed discount factor $v\in (0,1)$ and $x$ is the initial capital.

The description of the behavior of the surplus process is an amalgamation of those given in Subsections~\ref{sec:ci} and ~ \ref{sec:deF}, so we may be slightly more brief here. Namely, we stipulate that for $l\in \mathbb{N}$, during period $l$, a premium of $1$ is collected and the claim amount $\claim_l$ is incurred; then at the end of period $l$: (i) capital is injected in the amount $r_*(l)$, which is the amount by which the surplus process is negative ($r_*(l)=0$ if the surplus process remains nonnegative); (ii) the dividend amount $r(l)$ is paid out ($r(l)=0$ if the surplus process has become nonpositive). At end of period $0$ we inject $r_*(0)=(-x)\lor 0$ and a dividend $r(0)$ may be paid out, provided $x>0$. One says that the surplus process thus obtained is doubly reflected (at $0$ and $b$). The quantities paid out/injected at end of period $l$ are to be discounted by the factor $v^l$, $l\in \mathbb{N}_0$.

Then, using the fact proved in \cite{wu2011optimal} (\cite{APP} in the spectrally negative case), that barrier policies are optimal, the problem reduces to expressing, in terms of $W$ and $Z$, for a barrier dividend distribution policy $\pi_b$, $b\in \mathbb{N}_0$, the values of: (i) the expected discounted dividends,
\beq\label{e:Sobj}
  V^b_{D}(x):=\Ezb_x  \sum_{l=0}^\infty v^l r(l)=\bc
 \fr{Z_{v}(x)}  {\D Z_{v}(b)} & x \leq b\\
  x-b + V^b_{D}(b) & x > b\ec,  \eeq
  where   $\Ezb_x$ indicates expectation
with respect to the process doubly reflected at $0$ and $b$;
  and of (ii) the expected discounted bailouts,
    \beq \la{Bdisc} V^b_{B}(x):= \Ezb_x  \sum_{l=0}^\infty v^l r_*(l) =
\begin{cases}
Z_{v}(x)\fr{\D  Z_{1,v}(b)}  {\D Z_{v}(b)}  -Z_{1,v}(x) & x\leq b\\
V^b_{B}(b) &x>b
\end{cases}.
  \eeq
We give now the derivation of these two formulas. The cases $x>b$ are trivial, we limit the discussion to $x\leq b$.

\begin{proof}[Proof of \eqref{e:Sobj}, for dividends] 
We break the objective in two, following \cite{APP}: the ``De Finetti part" until   the first bailout time \eqref{e:Dobj}, and the rest \eqref{e:DP1}:
$$  V^b_{D}(x)=\Eb_x \left[ \sum_{l=0}^{\Tt-1} v^l r(l)\right]+  \Eb_x[ v^{\Tt};\Tt<\infty]  V^b_{D}(0)=\fr{W_v(x)}{\D W_v(b)}+ V^b_{D}(0)\le(Z_v(x) - \frac{\D Z_v(b)}{\D W_v(b)}W_v(x)\ri).$$
Making $x=0$ yields $V^b_{D}(0)=\frac 1{\D Z_v(b)}$ and the result follows.
 \end{proof}

\begin{proof}[Proof of \eqref{Bdisc}, for  bailouts] Using the strong Markov property at the exit time from the interval $[0,b)$ for the process $X$, yields an equation with three unknowns, $V^b_B(x)$, $V^b_B(0)$ and $V^b_B(b)$:
$$V^b_{B}(x)=\E_x  [v^{\tau_b^+};\tau_b^+ < \tau_{-1}^-]V^b_B(b)+\E_x  [v^{\tau^-_{-1}};\tau_{-1}^-<\tau_b^+]V^b_B(0) + \E_x[(-X(\tau_{-1}^-))v^{\tau_{-1}^-}; \tau_{-1}^-<\tau_b^+]$$
$$=\fr{W_{v}(x)}  {W_{v}(b)} V^b_{B}(b) + \le({Z_{v}(x)} -\fr{W_{v}(x)}  {W_{v}(b)} {Z_{v}(b)} \ri)V^b_{B}(0)-\le(Z_{1,v}(x) -\fr{W_{v}(x)}  {W_{v}(b)}  Z_{1,v}(b) \ri),$$
where, on the event that the first bailout occurs before the level $b$ is reached, the last term is the expectation  of this first bailout, before resetting to $0$, computed in Remark~\ref{r:deficit-at-ruin}, the penultimate term gives the expectation of the remaining bailouts and is given by \eqref{eq:Z-via-deficit-at-ruin}, finally the first term follows from \eqref{Wdef}.  Making $x=0$, yields $0= V^b_{B}(b) + Z_{1,v}(b) - {Z_{v}(b)} V^b_{B}(0)$, and it follows that
\footnotesize
$$   \fr{V^b_{B}(x)+ Z_{1,v}(x) - V^b_{B}(0) Z_{v}(x)}{W_{v}(x)}=\fr{V^b_{B}(b)+ Z_{1,v}(b) - V^b_{B}(0) Z_{v}(b)}{W_{v}(b)}=0.$$\normalsize
It remains to show that $V^b_{B}(0)=\fr{\D  Z_{1,v}(b)}  {\D Z_{v}(b)}$. To this end, using the strong Markov property at the exit time from the interval $[0,b]$ for the process $X$, produces $V^b_{B}(0)=\fr{W_{v}(0)}  {W_{v}(b+1)} V^b_{B}(b)+\le(1 -\fr{W_{v}(0)}  {W_{v}(b+1)} {Z_{v}(b+1)} \ri)V^b_{B}(0)+\fr{W_{v}(0)}  {W_{v}(b+1)}  Z_{1,v}(b+1)$. We conclude by plugging in $V^b_{B}(b) = {Z_{v}(b)} V^b_{B}(0)-Z_{1,v}(b)$.
\end{proof}
\beR
For $x<0$,  by Remark~\ref{r:deficit-at-ruin}, \eqref{Bdisc} reduces to $V^b_B(0)-x$, as it should.
\eeR
The combined  objective  is \begin{equation*} \la{cobj}    V(x):=\sup_{b\in \mathbb{N}_0} V^b(x),\quad V^b(x):= \Big[ V^b_{D}(x)- k V^b_{B}(x)\Big]=(x\lor b)-b+Z_v(x\land b)H(b)+k Z_{1,v}(x\land b),\end{equation*}
  with ``barrier influence function''
\begin{equation}
 H(b):=\fr{1- k \D  Z_{1,v}(b)}  {\D Z_{v}(b)}= \fr{1 -k \le(Z_{v}(b)-(1-\T p'(1-)) W_{v}(b)\ri)}
{(\frac{1}{v}-1) W_{v}(b)}. \la{e:H}
\end{equation}
As in the previous subsection, with an analogous justification, finding $V(x)$ is related to finding the supremum of $H$: (I) If $r:=\sup_{b\in \mathbb{N}_0}H(b)$ is attained, letting $b^*$ be a maximizer of $H$, then $b^*\geq x$ implies that $b^*$ is optimal for $V(x)$. (II) If the supremum defining $r$ is not attained, then the supremum defining $V(x)$ is not attained either. Since in this problem there is an optimal barrier strategy that does not depend on the initial reserve \cite[Theorem~3.2(B)]{wu2011optimal}, it follows, at least when the maximizer of $H$ is unique, that in (I), $b^*$ is in fact optimal for all $V(x)$, $x\in \Z$. Finally, note that $H$ differs from $\overline{H}(b):=\fr{1 -k Z_{v}(b)}{W_{v}(b)}$ only up to a positive affine transformation, so finding the supremum of, resp. a maximizer for, $H$ is equivalent to finding the supremum of, resp. a maximizer for, $\overline{H}$.

\section{
Examples}\label{sec:examples}

\ssec{Eventual ruin probabilities and the \deF dividends optimization}\la{s:edeF}
The  eventual ruin probability is a straightforward application of  \eqref{eq:perp-survival-transform}, followed by Taylor series coefficient extraction. Similarly, by using \eqref{eq:ruin-transform}-\eqref{disL} and generating function inversion, one can obtain the probability mass function of the time to ruin. One may also use the recursions  \eqref{dreci1}-\eqref{dreci2}-\eqref{eq:Lundberg}. Indeed, in the case when the support of the distribution of the claims is finite, the Lundberg recurrence \eqref{eq:Lundberg} reduces the problem of finding the eventual ruin probability to determining the roots of the characteristic equation \eqref{eq:characteristic-general}.

For instance, suppose $ \claim_1$ takes on the values $0$, $1$, and $3$, with probabilities $2/3$, $2/9$ and $1/9$, respectively. Then $E[ \claim_1]=5/9 < 1$ and eventual  upwards passage has probability
$\vf_1=1$. The \gf  is $\T p(z)=\fr{6}{9} + \fr{2}{9}
z + \fr{1}{9} z^3$, $z\in (0,1]$; 
 and Lundberg's equation is $$\fr{\vf_v}{v}=\fr{6}{9} + \fr{2}{9}
\vf_v + \fr{1}{9} \vf^ 3_v,\quad v \in (0,1].$$
The   recurrence for the perpetual survival and eventual ruin probabilities
writes as (with $f$ standing in place of $\Rui$ or $\sRui$)
$$f(x)= \fr{2}{3} f(x+1) +\fr{2}{9} f(x)+ \fr{1}{9} f(x-2)$$
for $x\in\mathbb{N}_0$. The characteristic equation \eqref{eq:characteristic-general} for this recurrence is (in $x$)
$$\fr{2}{3}x^3 -\fr{7}{9} x^2+ \fr{1}{9}=\fr{2}{3}\left(x-1\right)\left(x-\fr{1}{2} \right)\left(x+\fr{1}{3}\right)=0$$
(coinciding formally with the transformation of Lundberg's equation
$\fr{6}{9} - \fr{7}{9} \vf_1 + \fr{1}{9} \vf^ 3_1=0$,
 via $\vf_1\rightsquigarrow 1/x$).
 Satisfying the boundary conditions $\Rui(-1)=\Rui(-2)=1$, we arrive at $$\Rui(x)=\fr{2}{5} \left(\fr{1}{2}\right)^x-\fr{1}{15}
\left(-\fr{1}{3}\right)^x\text{ and }\sRui(x)=1-\fr{2}{5} \left(\fr{1}{2}\right)^x+\fr{1}{15}
\left(-\fr{1}{3}\right)^x,\quad x\in \mathbb{N}_0.$$ 
Taking $z$-transform yields, for $z\in (0,1)$,
$\T \Rui (z)=\frac{z+2}{(z+3)(2-z)}$ and $\T \sRui (z)= \frac{4}{6-7z+z^3}$,
which confirms \eqref{eq:z-transform-survival}-\eqref{PKd}. Finally, consider the \deF dividends optimization, under a
 discount factor  $v=150/169$.
Taylor expanding the scale transform \eqref{Wdisc} yields (the right-hand side features the consecutive values $\{W_v(0),W_v(1),\ldots\}$)\footnotesize
$$W_v=\{1.5, 2.035, 2.76082, 3.49551, 4.40307, 5.51337, 6.89721, 8.62338,
10.7802, 13.4755, 16.8446, 21.0558, 26.3198,
 ...\},$$ \normalsize which may be checked to be a convex  function with increasing  forward difference \footnotesize
$$\D W_v=\{0.535, 0.725817, 0.734691, 0.907565, 1.11029, 1.38385, 1.72616,
2.15678, 2.69539, 3.36905, 4.21121, 5.26398,  ...\}.$$\normalsize
   It follows that, irrespective of the initial capital, the optimal dividend policy is bringing the process to the barrier $b=0$ by taking dividends whenever possible.

\ssec{Modified geometric claims}\la{s:mgc}
We consider next modified
  geometric claims, defined by
$p_k= (1-\a) \a^{k-I} (1-p_0 -p_1-\cdots- p_{I-1}), \; k=I,I+1,...$, $\alpha\in [0,1)$. We restrict to $I=2$, which is equivalent to having two Lundberg roots \cite{sundt2007cramer}. We assume $p_0+p_1<1$. The \pgf is $$\T p(z)=\fr{p_0 + z(p_1 - \a p_0) + z^2 [(1-\a)(1-p_0)-p_1]}{1-\a z},\quad z\in (0,1].$$ The mean is  $m:=E \claim_1=1-p_0 +\fr{1-p_0-p_1}{1-\a}$, and the positive profit/subcritical case $m <1$ occurs when $ p_0(1-\a) >1-p_0-p_1$, which we assume henceforth. Fix $v\in (0,1]$. The Lundberg equation (in $z$) $k_vz^2 + z(p_1 - \a p_0- v^{-1})+ p_0=0$, with $k_v:=(1-\alpha)(1-  p_0) -p_1 +\alpha/v>0$,  has two (complex) solutions, the smaller one is $\vf_v$, and the larger of the two we will denote by $R_v$; their product is
$\vf_v R_v=\frac{p_0}{k_v}=\frac{p_0}{(1-\alpha)(1-  p_0) -p_1 +\alpha/v }$.

For $v=1$, the roots are $\vf_1=1$ and $R:=R_1=\fr{p_0}{1-p_1 -(1-\a)p_0}>1$. The eventual ruin probability is given by
\begin{equation}
\Rui(x)= \Rui(0) R^{-x}=\Big(1-\fr{1- E \claim_1}{p_0}\Big) R^{-x}=\fr{1- p_0 -p_1}{p_0(1-\a)} R^{-x}, \qu x\in \mathbb{N}_0. \end{equation}
This may be checked using \eqref{eq:z-transform-ruin}. Note the last formula does not hold for $x=-1$, except for special constellations of $p_0$, $p_1$, $\alpha$. Whatever the value of $v$, $\vf_v\leq 1<R_v$.

Some particular cases are:
\BEN \im If $\a=0$, the claims cannot exceed $2$, $R^{-1}=\fr{p_2}{p_0}$, $\Rui(0)= R^{-1}$, and the eventual ruin probability is
$$\Rui(x)=\left(\fr{p_2}{p_0}\right)^{x+1},\quad x\in \mathbb{N}_0\cup\{-1\},$$
recovering  the  classic gambler's ruin problem.
\im  Geometric: $p_0=1-\a$, $p_1=\a(1-\a)$.
\im  Geometric shifted by one: $p_0=0$, $p_1=1-\a$.
\im  Geometric shifted by two: $p_0=0= p_1$. \EEN

Writing now $\pg(z)-z/v=k_v(z-\f_v)(z-R_v)$, we find, using \eqref{Wdisc} \& \eqref{eq:Z-basic-transform}, for $z\in (0,\vf_v)$:
$$
\T W_v(z)=\frac{1}{k_v (R_v-\vf_v )}\Big(\frac{1}{\vf_v-z}-\frac{1}{R_v-z}\Big)$$ and
$$\T Z_v(z)=\frac{1}{1-z}\left(1-\frac{v^{-1}-1}{k_v(1-\vf_v)(R_v-1)}\right)+\frac{v^{-1}-1}{k_v(R_v-\vf_v)}\left(\frac{(\vf_v^{-1}-1)^{-1}}{\vf_v-z}+\frac{(1-R_v^{-1})^{-1}}{R_v-z}\right),
$$
so that, for $x\in \mathbb{N}_0$,
$$ W_v(x)=\frac{1}{k_v (R_v-\vf_v )}\Big(\vf_v^{-x -1}-R_v^{-x -1} \Big)$$ and
$$Z_v(x)=1-\frac{v^{-1}-1}{k_v(1-\vf_v)(R_v-1)}+\frac{v^{-1}-1}{k_v(R_v-\vf_v)}\left((\vf_v^{-1}-1)^{-1}\vf_v^{-x-1}+(1-R_v^{-1})^{-1}R_v^{-x-1}\right).$$
As a check, $\T W_v(0+)=W_v(0)=p_0^{-1}$ and $ \T Z_v(0+)=Z_v(0)=1$.
Given specific values of the parameters $\a$, $p_1$, $p_0$, $v$, the above expressions for $W_v$ and $Z_v$ may  be easily used to optimize combinations of expected bailouts/penalties at ruin and dividends.


\beR This model has a long history  in branching processes  as well \cite{AN,Mode}. Its utilisation there goes back to Steffensen and Lotka  (under the name of linear fractional branching) -- see \cite{kendall1966branching}, and  is still of interest  nowadays --  see for example \cite{Sag}.\eeR

\ssec {Multi-band dividend policies and modified \deF optimization}

 \beXa Recall 
 Morrill's historic example \cite[Ex. 2]{morrill1966one},
 with claims taking the values $0$ and $3$ with probabilities $12/13$ and $1/13$, respectively ($\rightarrow$ $E[\claim_1]=3/13<1$), and with discount factor  $v=65/72 $.
Taylor expanding the scale transform \eqref{Wdisc} yields \footnotesize
$$W_v=\{1.08333,1.3,1.56,1.78172,2.02973,2.30568,2.61834,2.97286,3.3753,3.83
   216,4.35085,...\},$$ \normalsize which may be checked to have a forward difference
\footnotesize $$\D W_v=\{0.216667,0.26,0.221722,0.248011,0.275947,0.312659,0.354523,0.402433,
   0.456864,...\}$$ \normalsize with a global minimum at $b^*= 0$
and another local minimum at  $2$.

\begin{figure}[h!]
\begin{center}
 \includegraphics[scale=.7]{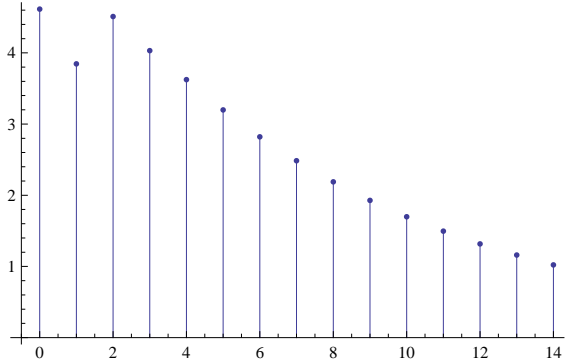}
\caption{The barrier influence function $1/\D W_v(b)$ for Morrill's example. The maximum $b^*=0$ is followed by the local maximum $2$. The optimal dividend policy  is multi-band, with two continuation sets $\{0\}$ and  $\{2\}$.
}\label{fig1}
\end{center}
\end{figure}

\begin{figure}[h!]
\begin{center}
 \includegraphics[scale=.7]{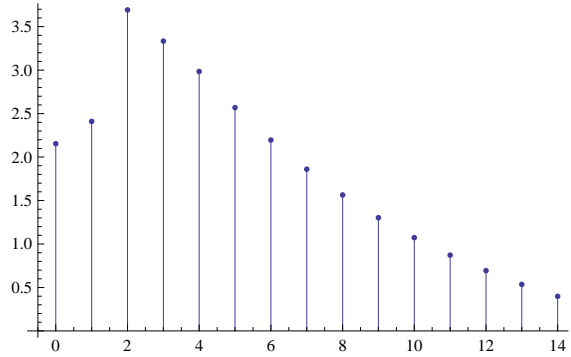}
\caption{The barrier influence function $(1-3.2 \D Z_{1,v}(b)) /\D W_v(b)$ for Morrill's example. The maximum $b^*=2$ is now the unique local (and hence global) maximum. 
}\label{fig2}
\end{center}
\end{figure}

Consider now the modified \deF objective of Subsection~\ref{sec:deficit-at-ruin-refl}. For $k$ big enough, for example   $k=3.2$, the  \bfun $(1-3.2 \D Z_{1,v}(b)) /\D W_v(b)$ is unimodal -- see Figure \ref{fig2}.

\eeXa

\begin{figure}[h!]
\begin{center}
 \includegraphics[scale=.7]{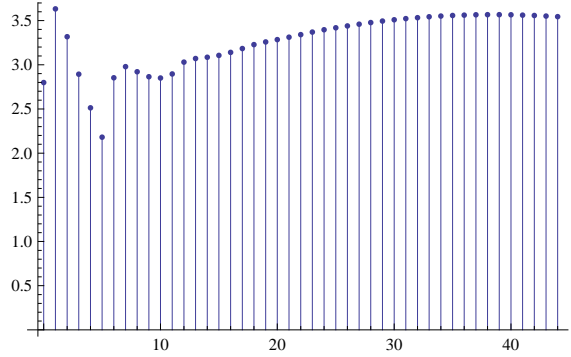}
\caption{The barrier influence function $1/\D W_v(b)$ for the Gerber-Shiu-Yang example. The maximum $b^*=1$ is followed by the local maxima $7$ and $38$. The optimal dividend policy  is multi-band.
}\label{fig3}
\end{center}
\end{figure}

\beXa\label{example:GSY} We turn now to the Gerber-Shiu-Yang example \cite[Ex.~3]{gerber2010elementary}, in which the barrier influence function $1 /\D W_v(b)$ has three local minima. The claims take the values $0$, $1$, $2$ and $7$ with probabilities $3/4,1/20,1/10$ and $1/10$, respectively ($\rightarrow$ $E[\claim_1]=19/20<1$), and the discount factor is $v=0.999 $. Now the \bfun has  a global maximum at $b^*= 1$ and two
further local maxima at  $ 7$ and $38$ --  see Figure \ref{fig3}.

\begin{figure}[h!]
\begin{center}
 \includegraphics[scale=.7]{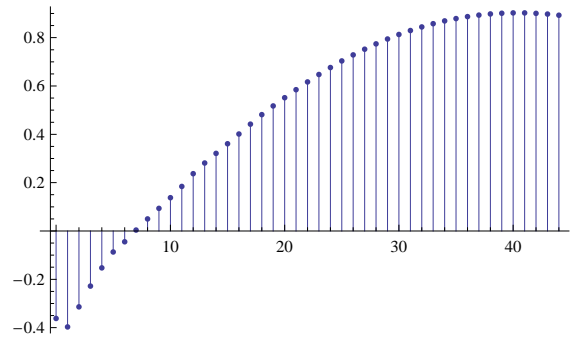}
\caption{The barrier influence function $(1-1.2 \D Z_{1,v}(b)) /\D W_v(b)$ for the Gerber-Shiu-Yang example. The unique local and global maximum is $b^*=41$.
}\label{fig4}
\end{center}
\end{figure}

\begin{figure}[h!]
\begin{center}
 \includegraphics[scale=.7]{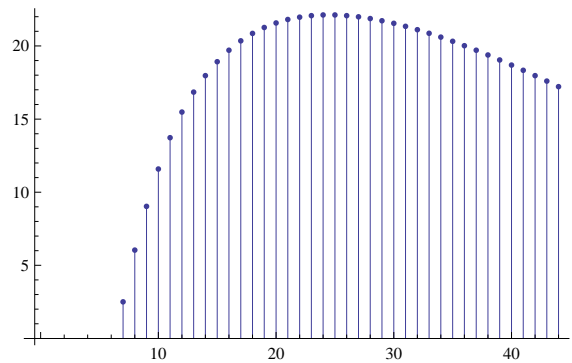}
\caption{The  barrier influence function \eqref{e:H} with $k=1.2$  for the doubly reflected Gerber-Shiu-Yang example has a unique global  maximum at $b^*=25$. The first seven  values $-89.91, -59.1845, -43.5339, -30.8171, -19.8565, -10.3512, -2.10264$ are too negative to be represented. The optimal policy is to take dividends above the level $b^*=25$. 
}\label{fig5}
\end{center}
\end{figure}

  Adopting a modified \deF objective of Subsection~\ref{sec:deficit-at-ruin-refl}, for example  with $k=1.2$ --- see  Figure~\ref{fig4}  --- shifts the global maximum to $b^*=41$.  The barrier influence function is not unimodal. 
With capital injections however, the \bfun is unimodal -- see next subsection and Figure~\ref{fig5}.

\eeXa
\ssec{Combined dividends and bailouts optimization objective for the doubly reflected process}
We optimize finally in the Gerber-Shiu-Yang example (Example~\ref{example:GSY}), the combined dividends-bailouts objective of Subsection~\ref{cobj} for the doubly reflected process with $k=1.2$. Recall that for this optimization problem there is always an optimal barrier policy \cite{wu2011optimal}. We obtain Figure~\ref{fig5}.
This objective seems to have achieved a ``compromise" between the peaks of the pure \deF objective.

\small
\bibliographystyle{alpha}
\bibliography{Pare34}

\begin{thebibliography}{AGVA17}

\bibitem[AGVA17]{AGV}
F.~Avram, D.~Grahovac, and C.~Vardar-Acar.
\newblock {The $W,Z$ scale functions kit for first passage problems of
  spectrally negative L\'evy processes, and applications to the optimization of
  dividends}.
\newblock {\em arXiv preprint arXiv:1706.06841}, 2017.

\bibitem[AI17]{AIjoint}
H.~Albrecher and J.~Ivanovs.
\newblock On the joint distribution of tax payments and capital injections for
  a {L\'e}vy risk model.
\newblock {\em Probability and Mathematical Statistics, to appear}, 2017.

\bibitem[AIZ16]{AIZ}
H.~Albrecher, J.~Ivanovs, and X.~Zhou.
\newblock Exit identities for {L\'e}vy processes observed at {P}oisson arrival
  times.
\newblock {\em Bernoulli}, 22(3):1364--1382, 2016.

\bibitem[AKP04]{AKP}
F.~Avram, A.~Kyprianou, and M.~Pistorius.
\newblock Exit problems for spectrally negative {L\'e}vy processes and
  applications to ({C}anadized) {R}ussian options.
\newblock {\em The Annals of Applied Probability}, 14(1):215--238, 2004.

\bibitem[AM05]{AM05}
P.~Azcue and N.~Muler.
\newblock Optimal reinsurance and dividend distribution policies in the
  {C}ram{\'e}r-{L}undberg model.
\newblock {\em Mathematical Finance}, 15(2):261--308, 2005.

\bibitem[AN72]{AN}
K.~B. Athreya and P.~E. Ney.
\newblock {\em Branching Processes}.
\newblock Springer, Berlin, 1972.

\bibitem[APP07]{APP}
F.~Avram, Z.~Palmowski, and M.~R. Pistorius.
\newblock On the optimal dividend problem for a spectrally negative {L\'e}vy
  process.
\newblock {\em The Annals of Applied Probability}, 17(1):156--180, 2007.

\bibitem[APP15]{APP15}
F.~Avram, Z.~Palmowski, and M.~R. Pistorius.
\newblock On {G}erber--{S}hiu functions and optimal dividend distribution for a
  {L\'e}vy risk process in the presence of a penalty function.
\newblock {\em The Annals of Applied Probability}, 25(4):1868--1935, 2015.

\bibitem[Ber97]{Ber97}
J.~Bertoin.
\newblock Exponential decay and ergodicity of completely asymmetric {L\'e}vy
  processes in a finite interval.
\newblock {\em The Annals of Applied Probability}, 7(1):156--169, 1997.

\bibitem[Ber98]{Ber}
J.~Bertoin.
\newblock {\em L{\'e}vy processes}.
\newblock Cambridge university press, 1998.

\bibitem[BF02]{banderier2002basic}
C.~Banderier and P.~Flajolet.
\newblock Basic analytic combinatorics of directed lattice paths.
\newblock {\em Theoretical Computer Science}, 281(1):37--80, 2002.

\bibitem[BJ15]{bauerle2015risk}
N.~B{\"a}uerle and A.~Ja{\'s}kiewicz.
\newblock Risk-sensitive dividend problems.
\newblock {\em European Journal of Operational Research}, 242(1):161--171,
  2015.

\bibitem[Bor12]{Bor}
A.~A. Borovkov.
\newblock {\em Stochastic Processes in Queueing Theory}.
\newblock Springer Science \& Business Media, 2012.

\bibitem[BPR10]{brown}
M.~Brown, E.~A. Pek{\"o}z, and S.~M. Ross.
\newblock Some results for skip-free random walk.
\newblock {\em Probability in the Engineering and Informational Sciences},
  24:491--507, 2010.

\bibitem[CGS00]{cheng2000discounted}
S.~Cheng, H.~U. Gerber, and E.~S.~W. Shiu.
\newblock Discounted probabilities and ruin theory in the compound binomial
  model.
\newblock {\em Insurance: Mathematics and Economics}, 26(2):239--250, 2000.

\bibitem[CKF06]{consul2006lagrangian}
P.~C. Consul, S.~Kotz, and F.~Famoye.
\newblock {\em Lagrangian Probability Distributions}.
\newblock Birkh{\"a}user Boston, 2006.

\bibitem[CP16]{PP}
M.~C.~H. Choi and P.~Patie.
\newblock Skip-free {M}arkov chains.
\newblock {\em Research gate}, 2016.

\bibitem[dF57]{DeF}
B.~de~Finetti.
\newblock Su un'impostazione alternativa della teoria collettiva del rischio.
\newblock In {\em Transactions of the XVth international congress of
  Actuaries}, volume~2, pages 433--443, 1957.

\bibitem[Don05]{Doney}
R.~A. Doney.
\newblock Some excursion calculations for spectrally one-sided {L\'e}vy
  processes.
\newblock In {\em S{\'e}minaire de Probabilit{\'e}s XXXVIII}, pages 5--15.
  Springer, 2005.

\bibitem[Fel71]{Fel}
W.~Feller.
\newblock {\em An Introduction to Probability Theory and its Applications},
  volume~II.
\newblock John Wiley \& Sons, New York, 1971.

\bibitem[Ger72]{gerber1972games}
H.~U. Gerber.
\newblock Games of economic survival with discrete-and continuous-income
  processes.
\newblock {\em Operations research}, 20(1):37--45, 1972.

\bibitem[Ger88]{gerber1988mathematical}
H.~U. Gerber.
\newblock Mathematical fun with ruin theory.
\newblock {\em Insurance: Mathematics and Economics}, 7(1):15--23, 1988.

\bibitem[GLY06]{gerber2006note}
H.~U. Gerber, X.~S. Lin, and H.~Yang.
\newblock A note on the dividends-penalty identity and the optimal dividend
  barrier.
\newblock {\em Astin Bulletin}, 36(02):489--503, 2006.

\bibitem[GSY10]{gerber2010elementary}
H.~U. Gerber, E.~S.~W. Shiu, and H.~Yang.
\newblock An elementary approach to discrete models of dividend strategies.
\newblock {\em Insurance: Mathematics and Economics}, 46(1):109--116, 2010.

\bibitem[IP12]{IP}
J.~Ivanovs and Z.~Palmowski.
\newblock Occupation densities in solving exit problems for {M}arkov additive
  processes and their reflections.
\newblock {\em Stochastic Processes and their Applications}, 122(9):3342--3360,
  2012.

\bibitem[Iva11]{Iva}
J.~Ivanovs.
\newblock {\em One-sided Markov additive processes and related exit problems}.
\newblock PhD thesis, Eurandom, 2011.

\bibitem[Kem61]{kemperman1961passage}
J.~H.~B. Kemperman.
\newblock {\em The Passage Problem for a Stationary Markov Chain}.
\newblock Statistical research monographs. University of Chicago Press, 1961.

\bibitem[Ken66]{kendall1966branching}
D.~G. Kendall.
\newblock Branching processes since 1873.
\newblock {\em Journal of the London Mathematical Society}, 1(1):385--406,
  1966.

\bibitem[KKR13]{KKR}
A.~Kuznetsov, A.~E. Kyprianou, and V.~Rivero.
\newblock The theory of scale functions for spectrally negative {L\'e}vy
  processes.
\newblock In {\em L{\'e}vy Matters II}, pages 97--186. Springer, 2013.

\bibitem[Kyp14]{Kyp}
A.~Kyprianou.
\newblock {\em Fluctuations of L{\'e}vy Processes with Applications:
  Introductory Lectures}.
\newblock Springer Science \& Business Media, 2014.

\bibitem[Loe08]{Loef}
R.~L. Loeffen.
\newblock {\em Stochastic control for spectrally negative L{\'e}vy processes}.
\newblock PhD Thesis, University of Bath, 2008.

\bibitem[Loe09]{loeffen2009optimal}
R.~L. Loeffen.
\newblock An optimal dividends problem with a terminal value for spectrally
  negative {L\'e}vy processes with a completely monotone jump density.
\newblock {\em Journal of Applied Probability}, 46(1):85--98, 2009.

\bibitem[LR10]{LR}
R.~L. Loeffen and J.-F. Renaud.
\newblock De {F}inetti's optimal dividends problem with an affine penalty
  function at ruin.
\newblock {\em Insurance: Mathematics and Economics}, 46(1):98--108, 2010.

\bibitem[Lun03]{Lun}
F.~Lundberg.
\newblock Approximerad framst{\"a}llning av sannolikhetsfunktionen.
\newblock {\em Akad. Afhandling. Almqvist och Wiksell}, 1903.

\bibitem[Mar01]{marchal2001combinatorial}
P.~Marchal.
\newblock A combinatorial approach to the two-sided exit problem for
  left-continuous random walks.
\newblock {\em Combinatorics, Probability and Computing}, 10(03):251--266,
  2001.

\bibitem[Miy61]{miyasawa1961economic}
K.~Miyasawa.
\newblock {\em An economic survival game}.
\newblock Econometric Research Program, Princeton University, 1961.

\bibitem[MM61]{miller1961dividend}
M.~H. Miller and F.~Modigliani.
\newblock Dividend policy, growth, and the valuation of shares.
\newblock {\em The Journal of Business}, 34(4):411--433, 1961.

\bibitem[Mod71]{Mode}
C.~J. Mode.
\newblock {\em Multitype branching processes: theory and applications},
  volume~34.
\newblock American Elsevier Publishing Company, 1971.

\bibitem[Mor66]{morrill1966one}
J.~E. Morrill.
\newblock One-person games of economic survival.
\newblock {\em Naval Research Logistics (NRL)}, 13(1):49--69, 1966.

\bibitem[MVJ15]{vidmar2015scales}
A.~Mijatovi{\'c}, M.~Vidmar, and S.~Jacka.
\newblock {Markov chain approximations to scale functions of L\'evy processes}.
\newblock {\em Stochastic Processes and their Applications},
  125(10):3932--3957, 2015.

\bibitem[Pis05]{Pispot}
M.~Pistorius.
\newblock A potential-theoretical review of some exit problems of spectrally
  negative {L\'e}vy processes.
\newblock {\em S{\'e}minaire de Probabilit{\'e}s XXXVIII}, pages 30--41, 2005.

\bibitem[Qui04]{quine}
M.~P. Quine.
\newblock On the escape probability for a left or right continuous random walk.
\newblock {\em {Annals of Combinatorics}}, 8:221--223, 2004.

\bibitem[Sag16]{Sag}
S.~Sagitov.
\newblock Tail generating functions for extendable branching processes.
\newblock {\em Stochastic Processes and their Applications}, 2016.

\bibitem[Sch07]{Schmidli}
H.~Schmidli.
\newblock {\em Stochastic control in insurance}.
\newblock Springer Science \& Business Media, 2007.

\bibitem[SdR07]{sundt2007cramer}
B.~Sundt and A.~D.~E. dos Reis.
\newblock Cram{\'e}r-{L}undberg results for the infinite time ruin probability
  in the compound binomial model.
\newblock {\em Bulletin of the Swiss Association of Actuaries}, 2:179--190,
  2007.

\bibitem[Shi89]{shiu1989probability}
E.~S.~W. Shiu.
\newblock The probability of eventual ruin in the compound binomial model.
\newblock {\em Astin Bulletin}, 19(2):179--190, 1989.

\bibitem[Spi13]{Spitzer}
F.~Spitzer.
\newblock {\em Principles of Random Walk}.
\newblock Springer Science \& Business Media, 2013.

\bibitem[Sup76]{Suprun}
V.~N. Suprun.
\newblock Problem of destruction and resolvent of a terminating process with
  independent increments.
\newblock {\em Ukrainian Mathematical Journal}, 28(1):39--51, 1976.

\bibitem[Tak77]{takacs1977combinatorial}
L.~Tak{\'a}cs.
\newblock {\em Combinatorial Methods in the Theory of Stochastic Processes}.
\newblock Wiley series in probability and mathematical statistics. R. E.
  Krieger Publishing Company, 1977.

\bibitem[Vid13]{vidmar2013fluctuation}
M.~Vidmar.
\newblock Fluctuation theory for upwards skip-free {L\'e}vy chains.
\newblock {\em arXiv preprint arXiv:1309.5328}, 2013.

\bibitem[Vid15]{vidmar2015overshoots}
M.~Vidmar.
\newblock {Non-random overshoots of L\'evy processes}.
\newblock {\em Markov processes and related fields}, 21(7):39--56, 2015.

\bibitem[WGT11]{wu2011optimal}
Y.~Wu, J.~Guo, and L.~Tang.
\newblock Optimal dividend strategies in discrete risk model with capital
  injections.
\newblock {\em Applied Stochastic Models in Business and Industry},
  27(5):557--566, 2011.

\bibitem[Wil93]{willmot1993ruin}
G.~E. Willmot.
\newblock Ruin probabilities in the compound binomial model.
\newblock {\em Insurance: Mathematics and Economics}, 12(2):133--142, 1993.

\bibitem[Xin04]{xin2004ring}
G.~Xin.
\newblock The ring of {M}alcev-{N}eumann series and the residue theorem.
\newblock {\em arXiv preprint arXiv:math/0405133}, 2004.

\end{thebibliography}
\normalsize
\appendix

\sec{Double (generating function) transforms of ruin probabilities \la{s:Wil}}
Recall the notation of Section~\ref{sec:non-smooth}. From \cite[Eqs.~(2.7) \& (2.13)]{willmot1993ruin}, one may deduced the double transform \beq  \la{PKdk}\T{\sRui}_v(z):=\sum_{n=0}^\I     v^n \sRui_z(n):=\sum_{n=0}^\I     v^n \le( \sum_{x=0}^\I z^x  \sRui(n;x)\ri)=\fr{\fr{z}{1-z}-\fr{\vf_v}{1-\vf_v} }{{z}-v {\pg(z)}}, \quad v, z \in (0,1),\; z\ne\vf_v,\eeq
where $\vf_v \in (0,1)$ is the Lundberg root \eqref{disL}
  (note that $z=\vf_v$  is a removable singularity). Indeed, from \eqref{dreci1}, for all $n\geq 1$, \cite[Eq.~(2.3)]{willmot1993ruin}
$$z\sRui_z(n)=\T p(z) \sRui_z(n-1) -p_0 \sRui(n-1;0),$$
and summing over $n$ after multiplication by $v^n$ yields  \cite[Eq.~(2.7)]{willmot1993ruin}
$$z(\T{\sRui}_v(z)-(1-z)^{-1})\!=\!v \T p(z)  \T{\sRui}_v(z)- p_0 v \sum_{n=0}^\I v^n \sRui(n;0) \Rightarrow (z-v \T p(z))\T{\sRui}_v(z)\!=\!\fr{z}{1-z}- p_0 v \sum_{n=0}^\I v^n \sRui(n;0),$$
from where \eqref{PKdk} is obtained by requiring that the root $z=\vf_v$  on the left-hand side annihilates also the right-hand side.

Eq.~\eqref{PKdk}
 implies the transform (for $v, z \in (0,1)$, $z\ne\vf_v$)

\begin{equation*}
\T{\Rui}_v(z):=\sum_{x=0}^\I \sum_{n=0}^\I  z^x v^n \Rui(n;x)=\fr 1{(1-z)(1-v)}-\T{\sRui}_v(z)=\fr{1 }{{z}-v {\pg(z)}}\le(\fr{v(z-\pg(z))}{(1-v)(1-z)}+\fr{\vf_v}{1-\vf_v}\ri).\end{equation*}

\beR
Note the single transforms \cite[Eq.~(3.5)]{willmot1993ruin}
\begin{eqnarray} \label{eq:z-transform-survival}
 \T{\sRui}(z)&:=&\sum_{x=0}^\I   z^x  \sRui(x)=\lim_{v\uparrow1} (1-v)\T{\sRui}_v(z)=\fr{(1-E [\claim_1])\lor 0}{\pg(z)-z},
\quad z\in (0,1),\\
\T{\Rui}(z)&:=&\sum_{x=0}^\I   z^x  \Rui(x)=\fr 1{1-z}-\fr{(1-E [\claim_1])\lor 0}{\pg(z)-z}, \la{PKd} \quad z\in (0,1), \label{eq:z-transform-ruin}\end{eqnarray}
which are similar to the \PK formulas of the \CL model. One also has \cite[2.14]{shiu1989probability} \begin{equation*}{\sRui}(0)=\lim_{z\downarrow  0}\T{\sRui}(z)=\fr{(1-E \claim_1)\lor 0}{p_0}.\end{equation*}
\eeR 
\sec{Summary table}\label{summary-table}
 \begin{center}
 \begin{tabular}{|p{4.8cm}|p{5.7cm}|p{5cm}|}
 \hline
\textbf{Right-continuous random walk}& \textbf{Upwards skip-free L\'evy chain} & \textbf{Spectrally negative L\'evy process} \\ \hline\hline
$\{n,b\}\subset \mathbb{N}_0$, $x\in \Z$, $x\leq b$, $v\in (0,1]$ & $\{q,\beta\}\subset [0,\infty)$, $\{b,x\}\subset h\Z$, $x\leq b$, $b\geq 0$ &  $\{q,\beta\}\subset [0,\infty)$, $\{b,x\}\subset  \mathbb{R}$, $x\leq b$, $b>0$\\\hline\hline
$X_n=X_0+n-\sum_{i=1}^nC_i$ & $Y=hX_{N}$, $N$ independent homogeneous Poisson process of intensity $\gamma$, $\{\gamma,h\}\subset (0,\infty)$ & L\'evy process $U$ having a.s. non-monotone paths and no positive jumps\\ \hline
$C_n$ i.i.d., $\mathbb{N}_0$-valued; p.m.f. $p$,
$p_0\in (0,1)$; p.g.f. $\pg$ & L\'evy measure $\lambda$; $\lambda=\gamma\sum_{i\in \mathbb{Z}\backslash \{1\}}p_i\delta_{h(1-i)}$; Laplace exponent $\psi$; $\psi(\beta)=\gamma[e^{\beta h}\pg(e^{-\beta h})-1]$& Laplace exponent $\psi$;  $\delta$ is the drift, when $X$ has bounded variation\\ \hline\hline
$\tau_b^-=\inf\{m \in \mathbb{N}_0: X_m\leq b\}$; & $\tau_b^-=\inf\{t\in [0,\infty): Y_t\leq b\}$; &$\tau_b^-=\inf\{t\in (0,\infty): U_t<b\}$; \\
$\tau_b^+=\inf\{m \in \mathbb{N}_0: X_t\geq b\}$&  $\tau_b^+=\inf\{t\in [0,\infty): Y_t\geq b\}$& $\tau_b^+=\inf\{t\in (0,\infty): U_t> b\}$\\\hline\hline
$\vf_v=$ smallest root of $\pg(\xi)/\xi=v^{-1}$ (in $\xi\in (0,1]$) & $\Phi(q)=$ largest root of $\psi(\lambda)-q$ (in $\lambda\in [0,\infty)$);  $e^{-h\Phi(q)}=\vf_{\frac{\gamma}{\gamma+q}}$ & $\Phi(q)=$ largest root of $\psi(\lambda)-q$ (in $\lambda\in [0,\infty)$)\\\hline
$ E_x \le[v^{\t_b^+}; \t_b^+<\infty\ri]=\vf_v^{b-x}$ &\multicolumn{2}{c|}{$E_x[e^{-q \tau_b^+};\tau_b^+<\infty]=e^{-\Phi(q)(b-x)}$} \\\hline\hline
$\sum_{y=0}^\I z^y W_v(y)=\fr {1}{\pg(z)-\frac{z}{v}}$, $z \in (0,\vf_v)$ & $\int_0^\infty e^{-\beta y}W^{(q)}(y)dy=\frac{e^{\beta h}-1}{\beta h(\psi(\beta)-q)}$, $\beta\in (\Phi(q),\infty)$;  $W^{(q)}$ c\`ad \& constant on each interval $[x,x+h)$; $W^{(q)}(x)=\frac{1}{\gamma h}W_{\frac{\gamma}{\gamma+q}}(x/h)$ & $\int_0^\infty e^{-\beta y}W^{(q)}(y)dy=\frac{1}{\psi(\beta)-q}$, $\beta\in (\Phi(q),\infty)$; $W^{(q)}$ continuous on $[0,\infty)$\\\hline
$E_x [v^{ \t_b};  \t_b^+ < \t_{-1}^-]=\fr{ W_v(x)}{ W_v(b)}$ & $E_x[e^{-q\tau_b^+};\tau_b^+<\tau_{-h}^-]=\frac{W^{(q)}(x)}{W^{(q)}(b)}$& $E_x[e^{-q\tau_b^+};\tau_b^+<\tau_0^-]=\frac{W^{(q)}(x)}{W^{(q)}(b)}$\\\hline
$ P_x(\tau_{-1}^-<\infty)=1-W_1(x)(1-\pg'(1-)\land 1) $& $P_x(\tau_{-h}^-<\infty)=1-W^{(0)}(x)(\psi'(0+)\lor 0)$& $P_x(\tau_0^-<\infty)=1-W^{(0)}(x)(\psi'(0+)\lor 0)$\\\hline
$W_v(x)=0$, $x<0$; $W_v(0)=1/p_0$ & $W^{(q)}(x)=0$, $x<0$; $W^{(q)}(0)=1/(h\lambda(\{h\}))$ & $W^{(q)}(x)=0$, $x<0$; $W^{(q)}(0)=0$ if $U$ has unbounded variation, $=1/\delta$ o/w\\\hline
$\lim_{y\to \infty}W_v(y)\vf_v^{y+1}=\frac{v}{1-v\pg'(\vf_v-)}$ & $\lim_{y\to\infty}W^{(q)}(y)e^{-\Phi(q)(y+h)}=1/\psi'(\Phi(q)+)$ & $\lim_{y\to\infty}W^{(q)}(y)e^{-\Phi(q)(y)}=1/\psi'(\Phi(q)+)$\\\hline\hline
$Z_v(x)=1+\left(\frac{1}{v}-1\right)\sum_{y=0}^{x-1}W_v(y)$, $x\geq 0$ & $Z^{(q)}(x)=1+q\sum_{k=0}^{x/h-1}W^{(q)}(kh)$, $x\geq 0$; $Z^{(q)}(x)=Z_{\frac{\gamma}{\gamma+q}}(x/h)$ & $Z^{(q)}(x)=1+q\int_0^xW^{(q)}(y)dy$, $x\geq 0$\\\hline
$Z_v(x)=1$, $x\leq 0$ & \multicolumn{2}{c|}{$Z^{(q)}(x)=1$, $x\leq 0$}\\\hline
$\sum_{y=0}^\I z^y  Z_v(y)=\fr{\T p(z)-z}{ (1-z)(\T p(z)-\frac{z}{v})}$, $z\in (0,\vf_v)$ & \multicolumn{2}{c|}{$\int_0^\infty Z^{(q)}(y)e^{-\beta y}dy=\frac{\psi(\beta)}{\beta(\psi(\beta)-q)}$, $\beta\in  (\Phi(q),\infty)$}\\\hline
$E_x [v^{ \t^-_{-1}}; \t^-_{-1} <\t_b^+]=Z_v(x)-\frac{W_v(x)}{W_v(b)}Z_v(b)$ & $E_x [e^{-q \t^-_{-h}}; \t^-_{-h} <\t_b^+]=Z^{(q)}(x)-\frac{W^{(q)}(x)}{W^{(q)}(b)}Z^{(q)}(b)$ & $E_x [e^{-q \t^-_0}; \t^-_0 <\t_b^+]=Z^{(q)}(x)-\frac{W^{(q)}(x)}{W^{(q)}(b)}Z^{(q)}(b)$\\\hline
$E_x[v^{\tau_{-1}^-};\tau_{-1}^-<\infty]=Z_v(x)-\frac{\vf_v(1-v)}{v(1-\vf_v)}W_v(x)$, $v<1$ & $E_x[e^{-q\tau_{-h}^-};\tau_{-h}^-<\infty]=Z^{(q)}(x)-\frac{qh}{e^{\Phi(q)h}-1}W^{(q)}(x)$, $q>0$ & $E_x[e^{-q\tau_0^-};\tau_0^-<\infty]=Z^{(q)}(x)-\frac{q}{\Phi(q)}W^{(q)}(x)$, $q>0$\\\hline\hline
$v^{m\land \tau_{-1}^m}W_v(X_{m\land \tau_{-1}^-})$ a martingale in $m\in \mathbb{N}_0$ & $e^{-q(t\land \tau_{-h}^-)}W^{(q)}(X_{t\land \tau_{-h}^-})$ a martingale in $t\in [0,\infty)$ & $e^{-q(t\land \tau_0^-\land  \tau_b^+)}W^{(q)}(X_{t\land \tau_0^-\land\tau_b^+})$ a martingale in $t\in [0,\infty)$\\\hline
 \end{tabular}
\end{center}

\beR
Every spectrally negative L\'evy process may be seen as a (weak) limit of a net $Y^h$ of upwards skip-free L\'evy chains, as $h\downarrow 0$ \cite{vidmar2015scales}. This means that a great many relations in the spectrally negative L\'evy setting may be got (at least naively) by simply passing to the limit $h\downarrow 0$ (formally, one must of course pay attention to whether or not the relevant functional is continuous with respect to such a weak limit).
\eeR

\beR
One of the important contributions of having a unified $\Phi,W,Z$ theory developed in all the three settings featuring in the table above, is that whenever a result is available for one of them, it may often be simply ``guessed'' in the others, by ``translating'' one set of quantities into the other (though ultimately it still needs to be proved). We have seen this time and again in the results of this paper. 
\eeR

\end{document}